\documentclass[aap,seceqn,MSNbibl,dvips]{arximspdf}

\doi{10.1214/10-AAP687}
\volume{20}
\issue{6}
\pubyear{2010}
\firstpage{2346}
\lastpage{2388}

\makeatletter
\newproclaim{rema}{Remark}
\newtheorem{lemma}{Lemma}

\newtheorem{prop}{Proposition}
\newtheorem{theorem}{Theorem}

\renewcommand{\l}{\lambda}

\def\sfrac#1#2{#1/#2}
\def\vfrac#1#2{(#1)/#2}
\def\afrac#1#2{#1/(#2)}
\def\vafrac#1#2{(#1)/(#2)}
\makeatother

\begin{document}
\begin{frontmatter}

\title{VRRW on complete-like graphs: Almost sure behavior}
\runtitle{VRRW on complete-like graphs}
\thankstext{tk}{Supported in part by NSERC research grant and by Alfred
P.~Sloan Research Fellowship.}
\thankstext{tr}{Supported by EPSRC Grant No.~EP/C01524X and European
Science Foundation.}
\begin{aug}
\author[A]{\fnms{Vlada} \snm{Limic}\corref{}\ead[label=e1]{vlada@cmi.univ-mrs.fr}} and
\author[B]{\fnms{Stanislav} \snm{Volkov}\ead[label=e2]{S.Volkov@bristol.ac.uk}}
\runauthor{V. Limic and S. Volkov}
\affiliation{CNRS and Universit\'e de Provence and University of Bristol}
\address[A]{CNRS and Universit\'e de Provence\\
UMR-CNRS 6632, C.M.I.\\
13453 Marseille Cedex 13\\
France\\
\printead{e1}} 
\address[B]{Department of Mathematics\\
University of Bristol\\
BS8~1TW Bristol\\
UK\\
\printead{e2}}
\end{aug}

\received{\smonth{4} \syear{2009}}
\revised{\smonth{1} \syear{2010}}

%
\begin{abstract}
By a theorem of Volkov \cite{V2001} we know that on most graphs with
positive probability the linearly vertex-reinforced random walk\break (VRRW)
stays within a finite ``trapping'' subgraph at all large times.
The question of whether this tail behavior occurs with probability
one is open in general. In his thesis, Pemantle \cite{Pemthesis} proved,
via a dynamical system approach,
that for a VRRW on any complete graph the asymptotic frequency of visits
is uniform over vertices. These techniques do not easily extend even to the
setting of complete-like graphs, that is, complete graphs
ornamented with finitely many leaves at each vertex. In this work
we combine martingale and large deviation techniques to prove that
almost surely the VRRW on any such graph spends positive (and equal)
proportions of time on each of its nonleaf vertices. This behavior was
previously shown to occur only up to event of positive
probability (cf.~Volkov \cite{V2001}). We believe that our approach can
be used as a building block in studying related questions on more
general graphs. The same set of techniques is used to obtain explicit
bounds on
the speed of convergence of the empirical occupation measure.
\end{abstract}

%

\end{frontmatter}

\section{Introduction}\label{Introduction}

Consider a complete-like graph $\mathcal{G}_d$ with $d\geq2$ \textit{interior vertices} (or sites) and $r_i\geq0$ exterior vertices or
\textit{leaves} attached to the $i$th interior site, $i\in
\{1,\ldots,d\}$. More precisely, denote by
$V_d=\{1,2,\ldots,d,\ell_1^1,\ldots,\ell_{r_1}^1,\break\ldots,\ell
_1^d,\ldots,\ell_{r_d}^d\}$
the set of sites of $\mathcal{G}_d$, and by $E_d$ the set of its edges.
Typically we denote the edge connecting two different sites $v$
and $w$ by $\{v,w\}$. Any two sites that share an edge are called
\textit{neighbors}. If $v$ and $w$ are neighbors we also write
$v \sim w$. Then $E_d$ consist of $d(d-1)/2$ edges connecting
each pair of interior sites, as well as of the edges
$\{i,\ell_r^i\}$, for each $i\in\{1,\ldots,d\}$ and
$r=1,\ldots,r_i$. We will refer to $\ell_r^i$ as the $r$th leaf
attached to the interior vertex $i$. It is possible that $r_i=0$
for some $i$, in which case there is no leaf attached to $i$. If
$r_i=0$, for all $i=1,\ldots,d$, then $\mathcal{G}_d$ is the complete
graph on $d$ vertices. Any graph from the above class can be
viewed as a ``perturbation'' of the complete graph.

We start by recalling the (discrete-time) linearly vertex reinforced random
walk (VRRW) (see, e.g.,~\cite{Review}). This process can be
constructed on
general bounded degree graphs, but since the current work concerns
VRRW on complete-like graphs given above, the definition below can
be read with this special setting in mind.

The \textit{time} $t$ will run through positive integers. We denote
by $X(t)$ the position (site) of the walk at time $t$. Assume that
$z(0,v)$ are given positive integer quantities. For example, it
could be $z(0,v)\equiv1$, $v\in V_d$. Without loss of generality,
we can assume that the initial time is $t_0 = \sum_{v\in V_d}
z(0,v)$. Let $Z(t,v)$ equal $z(0,v)$ plus the number of visits
to vertex $v\in V_d$ up to time $t$, $t\geq t_0$.
Note that in this way we have $\sum_{v\in V_d} Z(t,v)\equiv t$ for
$t\geq t_0$. Denote by $(\mathcal{F}_t,t\geq t_0)$ the filtration generated
by $(X(t),t\geq t_0)$ (or equivalently by $(Z(t,v),t\geq t_0)$,
$v\in V_d$) up to time $t$. Then on the event $\{X(t)=v\}$ the
transitions of our process are given by
\begin{eqnarray}\label{Etrans}
\mathbb{P}\bigl(X(t+1) =w |\mathcal{F}_t\bigr)= \frac{Z(t,w)}{\sum_{y\in
V_d:y \sim v}
Z(t,y)}
\end{eqnarray}
for all $w\in V_d$, $w \sim v$. In particular, when at $\ell_r^i$,
the walk must return to $i$ in the next step.

Let
\begin{eqnarray*}
\pi(t)&=&\frac{1}{t}(Z(t,1),Z(t,2),\ldots,Z(t,d),
\\
&&\hspace*{12pt}{}Z(t,\ell_1^1),\ldots,Z(t,\ell_{r_1}^1),\ldots, Z(t,\ell_1^d),\ldots,Z(t,\ell_{r_d}^d))
\end{eqnarray*}
be the occupation measure generated by the VRRW above at time $t$,
determined by the vector of its atoms.
Let
$\pi_{\infty}=\lim_{t\to\infty} \pi(t)$ be the
asymptotic occupation measure on the event where
this limit exists, and
set $\pi_{\infty}=(0,0,\ldots,0)$ on the complement. Note that
$\pi(t)\in{\mathbb{R}}^{|V_d|}$, for all $t$, where
$|V_d|:=d+\sum_{i=1}^{d}r_i$, and we use this fact without further
mention. Set
\[
\pi_{\mathsf{unif}}:= \biggl(\frac1d,\frac1d,\ldots,\frac1d, 0, \ldots
,0 \biggr),
\]
where the initial $d$ coordinates are positive, and the
other $\sum_{i=1}^{d}r_i$ are equal to $0$.

The first goal of this paper is to prove
\begin{theorem}\label{th1}
For VRRW on $\mathcal{G}_d$, $d\ge3$, we have $\mathbb{P}(\pi
_{\infty}=\pi
_{\mathsf{unif}})=1$.
\end{theorem}

The next statement is related to the slow speed of convergence
noticed by Pemantle and Skyrms in \cite{PeSk}. Denote by $\Vert
\cdot\Vert=\Vert\cdot\Vert_{\infty}$ the maximum norm on
${\mathbb{R}}^{|V_d|}$.
\begin{theorem}\label{th2}
Let $\mathcal{G}_d$ be the complete-like graph on $d\geq3$ vertices. Then
for any $\delta>0$
\begin{eqnarray}\label{Eth2uppera}
\mathbb{P}\Bigl(\limsup_{t\to\infty} \Vert\pi(t)-\pi_{\mathsf{unif}}\Vert t^{\sfrac13-\delta
}<\infty\Bigr)&=&1 \qquad\mbox{if }d=3, 4,\\
\label{Eth2upperb}
\mathbb{P}\Bigl(\limsup_{t\to\infty} \Vert\pi(t)-\pi_{\mathsf{unif}}\Vert t^{\afrac 1{d-1}}<\infty\Bigr)&=&1 \qquad\mbox{if } d\ge5.
\end{eqnarray}
Moreover, for each $d\geq3$, if $|V_d|\geq d+1$ (there exists at least
one leaf) and any $\delta>0$
\begin{eqnarray}\label{Eth2lower}
\mathbb{P}\Bigl(\liminf_{t\to\infty} \Vert\pi(t)-\pi_{\mathsf{unif}}\Vert t^{\vafrac{d-2}{d-1}+\delta}=\infty\Bigr)=1.
\end{eqnarray}
\end{theorem}

In particular, the empirical occupation measure converges to
$\pi_{\mathsf{unif}}$ at least as fast as an inverse of a certain
power function,
and not faster than an inverse of another power function (provided
$|V_d|>0$). Note that (\ref{Eth2lower}) gives an upper bound on
the power exponent which is strictly smaller than $1$. To the best
of our knowledge, this is the first rigorous result verifying
``slow convergence'' for this class of models. However, the
problem of finding a lower bound on the speed in the case of the
complete graph is still open, and we believe that the true rate of
convergence is closer to the one in (\ref{Eth2uppera}) and (\ref
{Eth2upperb}). We wish to point out that computer simulations seem
to be misleading in predicting/confirming any of the above
results, due to the slow speed of convergence. With this in mind,
it is worth mentioning that our computer simulations seem to
suggest that for $d=3$
\[
\frac{\log{\mathsf M} (\Vert\pi(t)-\pi_{\mathsf{unif}}\Vert)}{\log t}\to-\frac12,
\]
where ${\mathsf M}(X)$ stands for the median of a random variable $X$.
The special case $d=2$ will be discussed in Section
\ref{S:casetwo}.

There exist a few mathematical results on the
asymptotic behavior of VRRW preceding this work.
As mentioned in the abstract,
Pemantle \cite{Pemthesis} proved that on any complete graph the
asymptotic frequencies of visits by the VRRW are the same for all
vertices. The papers \cite{PV1999} and \cite{Tar2004} study the
VRRW on the integers $\mathbb{Z}$. Pemantle and Volkov \cite{PV1999} prove
that this VRRW cannot get trapped on a subgraph spanned by $4$
sites, and moreover that it gets trapped on a random subgraph
spanned by $5$ subsequent sites with a positive probability.
Tarr\`es \cite{Tar2004} proved that this striking behavior occurs
almost surely, using subtle martingale and coupling techniques.

A study by Volkov \cite{V2001} exhibits a family of ``trapping
subgraphs'' for the VRRW on a general graph, where the range of
the VRRW is contained in any such subgraph. Recent results of
Bena\"im and Tarr\`es \cite{BenTar} show similar localization
phenomenon for certain natural generalizations of VRRW. The
asymptotic results in both \cite{BenTar} and \cite{V2001} are
shown to hold only on an event of \textit{positive} probability.
Volkov \cite{V2006} initiated the analysis of nonlinearly
reinforced VRRW. His analysis mostly concentrated on the power-law
reinforcement functions and the VRRW on $\mathbb{Z}$. Many interesting
open questions remain.

The rest of the paper is organized as follows.
Sections~\ref{S:polyaexample}--\ref{securn} recall a few
techniques used in related settings, and establish some
preliminary results. In Section~\ref{sec_MT} we introduce a
modified VRRW on a triangle with one special (more reinforced)
vertex and study the asymptotics of weights on the nonspecial
vertices. Section~\ref{S:anaclike} contains the proof of Theorem
\ref{th1} in the general (and novel) case of complete-like graphs
$\mathcal{G}_d$, and Section~\ref{Sec:generalize} discusses some
generalizations for $d$-partite graphs with leaves. Finally, in
Section~\ref{S:speed} we show Theorem \ref{th2}.

We will use the symbol $\wedge$ (resp.,~$\vee$) to denote the
operation of taking the minimum (resp.,~maximum) of two or more
numbers. For $f$ and $g$, two sequences of positive functions
defined on the positive reals, we write $f(t)=O(g(t))$ if
$\limsup_t f(t)/g(t)$ is finite, $g(t) \asymp f(t)$ or $f(t)=\Theta
(g(t))$ if both
$f(t)=O(g(t))$ and $g(t)=O(f(t))$, and $f(t)=o(g(t))$ if $\lim_t
f(t)/g(t)=0$.
The above notations extend in a straightforward way
to the stochastic setting.

\subsection{Multi-color {P\'{o}lya} urns and VRRW on complete graphs}\label{S:polyaexample}

We devote this short subsection to a
calculation that will hopefully both stimulate the reader's
interest in the problem, and point out some of the difficulties
awaiting. In addition, we will use a modification of the
supermartingale below in arguments of Section~\ref{S:anaclike}.
Fix $d\geq2$, and let $\Pi$ be the $d$-color {P\'{o}lya} urn started
with one ball of each color. In particular, at each step, one ball
is drawn from the urn at random, and it is placed back immediately
together with another ball of the same color. As usual, let the
initial time be $d$, and for each time $t\geq d$ denote by
$\Pi_i(t)$ the number of balls of color $i$, $i=1,\ldots,d$ in the
urn at time $t$. In this way $\sum_{i=1}^d \Pi_i(t) = t$ always. A
slick way (see \cite{V2001}, Section 2.1) to prove convergence of
the frequencies $\Pi_i(t)/t$, $i=1,\ldots,d$, to nontrivial
(nonzero, a.s.) random variables is via the following martingale
method. Using classical martingales $\Pi_i(t)/t$ for showing this
convergence is not optimal for showing that the limit is nonzero,
almost surely. Define
\[
M_i(t):=\log(t)-\log\bigl(\Pi_i(t)-1\bigr),
\]
and then check that the drift of this process equals
\[
\mathbb{E}\bigl(M_i(t+1)-M_i(t)| \mathcal{F}_t\bigr)
= \log\biggl(1+\frac{1}{t} \biggr)-\frac{\Pi_i(t)}{t}\log\biggl(1+\frac{1}{\Pi_i(t)-1} \biggr),
\]
and is therefore almost surely negative. Thus $M_i(t)$ is a
nonnegative supermartingale and it converges almost surely to a
finite quantity, hence $\Pi_i(t)/t$ converges almost surely to
a positive quantity.

Next consider the VRRW on complete graph with $d$ vertices. The
only difference of transitions of $(Z(t,1),\ldots,Z(t,d))$ from
those of $(\Pi_1(t),\ldots,\Pi_d(t))$ is that $\Pi_i(t+1)$ becomes
$1+\Pi(t)$ with probability proportional to $\Pi_i(t)$ no matter
which ball was drawn at time $t-1$, while $Z(t+1,i)$ becomes
$1+Z(t,i)$ with probability proportional to $Z(t,i)$ only if the
current position of the VRRW is not $i$; in turn this proportion
is taken with respect to the values at all but the currently
visited site. If one tries simply to recycle the above
supermartingale by subtracting a drift increment of order $1/t$ at
each time $t$ when $Z(t,i)=Z(t-1,i)+1$, then on the event that
$Z(t,i)$ is asymptotically of order larger than $t/\log(t)$
[this happens, since $Z(t,i) \sim t/d$, a.s.] the sum of the drift
increments diverges and it not possible to conclude convergence of
$M_i(t)$. One could think that there should be a simple way to
overcome the above difficulty, but we are not aware of one.

\subsection{Large deviation tools}\label{secLD}

Part of our analysis
(cf.~Section \ref{secproofs}) will use the strategy of
Volkov \cite{V2001} (see also 
\cite{BenTar}).

We recall the following classical facts.
Let $\xi_i$ be i.i.d. random variables with
$\mathbb{P}\{\xi_i=1\}=1-\mathbb{P}\{\xi_i=0\}=p\in(0,1)$. Define
for $a,p\in
(0,1)$,
%
\begin{equation}\label{entro}
H(a,p):=a\log\frac ap+(1-a)\log\frac{1-a}{1-p}\geq0.
\end{equation}
Recall an elementary fact from large deviation theory (see, e.g.,
\cite{SHIR}): for any $a^+\in[p,1)$ and any $a^-\in(0,p]$, we
have
%
\begin{equation}\label{LDPg}
\hspace*{20pt}\mathbb{P}\Biggl\{\frac1n \sum_{i=1}^n \xi_i\geq a^+ \Biggr\}\leq e^{-nH(a^+,p)},\qquad
\mathbb{P}\Biggl\{\frac1n \sum_{i=1}^n \xi_i\leq a^- \Biggr\}\leq e^{-nH(a^-,p)}.
\end{equation}

It is easy to verify (see also Propositions 2.2 and 2.3 in \cite{V2001}) that
%
\begin{eqnarray}\label{Hpro1}
 H(a,p)&=&\frac{\delta^2}{2p(1-p)}+\Theta\biggl(\frac{\delta^3}{p^2(1-p)^2} \biggr) \qquad\nonumber
\\
\eqntext{\mbox{if }a=p \pm\delta,\mbox{ where }\delta\ll1\quad \mbox{and}\nonumber} \\[-8pt]\\[-8pt]
 H(a,p)&=&p(r\log r -r +1)+\Theta(p^2 )\nonumber
\\
\eqntext{\mbox{if }a=rp,r=\Theta(1),\mbox{ and }a\vee p \ll1.\nonumber}
\end{eqnarray}

\subsection{Urn and martingale tools}\label{securn}

We start by recalling the results on urns from Pemantle and Volkov
\cite{PV1999}. We will often use them directly in coupling
arguments; however we will also need to generalize Theorem \ref{thurn1}
below (see Lemma~\ref{Ljointexc}) during the course of our
analysis.

The urn model defined below generalizes both the (original) {P\'{o}lya}
and the Friedman urn, and it is sometimes referred to as the \textit{generalized {P\'{o}lya} urn}. Consider the dynamics
\begin{eqnarray}\label{equrn}
(X_{n+1}, Y_{n+1}) &=& (X_n + a, Y_n + b)\qquad \mbox{with probability }\frac{X_n}{X_n + Y_n},\nonumber \\[-8pt]\\[-8pt]
(X_{n+1}, Y_{n+1}) &=& (X_n + c, Y_n + d)\qquad \mbox{with probability } \frac{Y_n}{X_n + Y_n}. \nonumber
\end{eqnarray}
We do not necessarily assume that the random numbers
$X_n,Y_n$ (of balls) are integer valued. When $\bigl({a \atop c} \enskip{b \atop d}\bigr)$ is a
multiple of the identity matrix (resp.,~$a=d$ and $b=c$ are all
nonzero), we recover {P\'{o}lya}'s (resp.,~Friedman's) urn. In all
cases where $\bigl({a \atop c} \enskip{b \atop d}\bigr)$ has a left
eigenvector $(v_1, v_2)$ with
positive components and $abcd > 0$, Freedman's
analysis~\cite{Freedman} can be carried through to show that $X_n
/ (X_n + Y_n)$ converges a.s. to $v_1 / (v_1 + v_2)$. When $a
> d, b > 0$ and $c = 0$ the urn is still Friedman like: although
$(0,1)$ is an eigenvector, it is easy to see that the principal
eigenvector is $(a-d, b)$ and that $X_n / (X_n + Y_n) \rightarrow
(a-d) / [(a-d) + b]$ a.s. The case $ad = bc = 0$ is trivial, so we
are left with the cases $ad > 0 = b = c$ and $ad
> 0 = bc < b + c$. Multiplication of $\bigl({a \atop c} \enskip{b \atop d}\bigr)$ by a constant does
not affect the asymptotic behavior. Due to symmetry, the
interesting behavior is captured in the following two theorems.
\begin{theorem}[(\cite{PV1999}, Theorem 2.2)]\label{thurn1}
Suppose $a > d = 1$, and $b = c = 0$. Then $\log X_n / \log Y_n
\rightarrow a$.
\end{theorem}
\begin{theorem}[(\cite{PV1999}, Theorem 2.3)]\label{thurn2}
Suppose $a = d = 1$, $b = 0$ and $c > 0$. Then $X_n / (c Y_n) - \log
Y_n$ converges to a random limit in $(-\infty, \infty)$.
\end{theorem}

\begin{rema}
(1) Theorem~\ref{thurn1} implies that for any $\varepsilon>0$ we have
$X_n^{(1/a-\varepsilon)} \leq Y_n \leq X_n^{(1/a+\varepsilon)}$ for
all large
$n$, almost surely. Since $X_n+Y_n \asymp n$, this easily implies
that $X_n$ is equal to $a\cdot n$ plus lower order terms, while
$Y_n$ is asymptotically equal to $n^{1/a}$ multiplied by a random
factor $A_n$, where for any $\varepsilon>0$ $A_n\in(n^{-\varepsilon
},n^\varepsilon)$
for all
large $n$.

(2)
The result in Theorem~\ref{thurn2} may be more surprising, in
that it shows $Y_n$ to be of the order $n / \log n$ multiplied by
a specific constant, with a random lower order correction. That
is, $X_n$ is asymptotically $c Y_n (A + \log Y_n)$, where $A$ is a
random constant. This class of urns was used in \cite{PV1999} to
prove that VRRW on $\mathbb{Z}$ cannot get trapped on a subgraph spanned
by $4$ subsequent points. Note that in the special case $c = 1$,
the urn process corresponds to a VRRW on the graph $\mathcal{G}$ with
$V(\mathcal{G}) = \{ u, v \}$, having one edge between $u$ and $v$ and
one loop connecting $u$ to itself, observed at the times of
successive visits to vertex $u$. Thus VRRW on this $\mathcal{G}$ spends
roughly
$n / \log n$ units of time at $v$ up to time $n$.

(3) Both of the above theorems can be derived using an elegant
method of Athreya and Ney \cite{ANey}, by embedding the urn into
a continuous time multi-type branching process. However, the proof
by embedding (see also  \cite{Janson} for recent progress)
is much less robust to ``variations'' in dynamics than the
martingale proofs of \cite{PV1999}. One such
variation is the setting where some (or all) of the parameters
$a,b,c,d$ are perturbed about fixed values (their means), and
where the distribution of these random perturbations varies over
time. Section~\ref{sec_MT} is devoted to proving some extensions
in this direction that turn out to be essential for our analysis.
\end{rema}

In the current work, we will repeatedly
bound the $\limsup$ (by a finite random quantity) of a process
that has supermartingale increments whenever its value is sufficiently large
via a separate martingale technique (see
Chapter 4 of \cite{TarThesis} for a similar idea in
a somewhat simpler setting).

In our general setting,
we are given $(\xi_n, n\geq0)$, a discrete-time process (not
necessarily bounded below nor above),
adapted to a filtration $(\mathcal{F}_n, n\geq0)$. In addition,
suppose there
exists $a,b\in{\mathbb{R}}$, $b>0$ such that:
\begin{enumerate}[(1)]
\item[(1)]
$ \xi$ has supermartingale increments on $[a,\infty)$, that is,
%
\begin{equation}
\label{Efirstp}
\mathbb{E}\bigl((\xi_{k+1}-\xi_k) 1_{\{\xi_k\geq a\}}|\mathcal{F}_k\bigr)
\leq0;
\end{equation}
\item[(2)] the overshoot of $\xi$ across $a$ is asymptotically
bounded by $b$, that is,
%
\begin{equation}
\label{Esecondp}
o^*(a):=\limsup_k
1_{\{\xi_k<a<\xi_{k+1}\}}
(\xi_{k+1} -a) \leq b
\qquad\mbox{almost surely};
\end{equation}
\item[(3)]
the tail variance of $\xi$ on $[a,\infty)$ is finite, that is,
%
\begin{equation}
\label{Ethirdp}
\sum_{k} \mathbb{E}\bigl[(\Delta\xi_k)^2 1_{\{\xi_k\geq a\}}\bigr]<\infty
\qquad\mbox{where }\Delta\xi_k:=\xi_{k+1}-\xi_k.
\end{equation}
\end{enumerate}

\begin{lemma}\label{Ljointexc}
Under the above assumptions
\[
\xi^*:=\limsup_{n\to\infty} \xi_n <\infty,\qquad \mbox{a.s.}
\]
\end{lemma}

\begin{pf}
Due to shift and scaling, without loss of generality
(WLOG) we may assume that $a=-1$ and $b=1$. Next fix a small $\delta
>0$, and define
\[
B_\delta^{(n)}=\Bigl\{\sup_{k\geq n} 1_{\{\xi_k<-1<\xi_{k+1}\}}
\bigl(\xi_{k+1}-(-1)\bigr) \leq1+\delta\Bigr\}.
\]
Property (\ref{Esecondp}) can be restated as $\lim_{n\to
\infty}\mathbb{P}(B_\delta^{(n)})= 1$. We shall now introduce an auxiliary
process
\[
{\xi'}^{,(n,\delta)}\equiv xi' :=(\xi'_k, k\geq n),
\]
adapted to the filtration generated by $(\xi_k, k\geq n)$,
and such that the three properties (\ref{Efirstp})--(\ref{Ethirdp})
hold for
$\xi'$, with $a=\delta$ and $b=0$.
Moreover, the inequality in (\ref{Efirstp}) for $\xi'$ becomes equality
%
\begin{equation}
\label{Efirstpeq}
\mathbb{E}\bigl((\xi_{k+1}'-\xi_k') 1_{\{\xi_k'\geq\delta\}
}|\mathcal{F}_k\bigr) = 0,\qquad k \geq n,
\end{equation}
and also
%
\begin{equation}
\label{Eonedominother}
B_\delta^{(n)}\subset\bigcap_{k\geq n} \{\xi_k\leq\xi_k'\}
\qquad\mbox{almost surely}.
\end{equation}
Define $\xi'_n\equiv{\xi'}^{,(n,\delta)}_n:=\xi_n$, and for
$k\geq n$ let
%
\begin{equation}
\label{Ecoupling}
\xi_{k+1}' := \cases{
\xi'_k+ \Delta\xi_k - \mathbb{E}(\Delta\xi_k|\mathcal{F}_k),
&\quad\mbox{if }$\xi_k \geq-1$,\cr
(\xi'_k+ \Delta\xi_k)\wedge\delta, &\quad\mbox{if }$\xi_k<-1$\mbox{ and }$\xi_k'< \delta$,\cr
\xi'_k, &\quad\mbox{if }$\xi_k <-1$\mbox{ and }$\xi_k'\geq\delta$.
}
\end{equation}
If $\xi_k'\geq\delta$ then either $\xi_k\geq-1$ in which case the
increment of $\xi'$ is the Doob--Meyer martingale ``correction'' of
the increment of~$\xi$, or $\xi_k< -1$ and then $\xi'$ does not
change value. So indeed, (\ref{Efirstp}) holds for $\xi'$ as
(\ref{Efirstpeq}). The property (\ref{Esecondp}) is immediate
since a positive overshoot of $\xi'$ across $\delta$ may occur only as
a result of a jump of $\xi$ when its current value is greater than
$-1$, but these jumps are asymptotically negligible by
(\ref{Ethirdp}). Similarly, (\ref{Ethirdp}) for $\xi'$ is easy to
derive from the definition (\ref{Ecoupling}), the property
(\ref{Ethirdp}) for\vspace*{1pt} $\xi$, and the standard fact ${\mathbb
{E}}((\Delta\xi_k
- \mathbb{E}(\Delta\xi_k|\mathcal{F}_k))^2 |\mathcal{F}_k)\leq
\mathbb{E}((\Delta\xi_k)^2
|\mathcal{F}_k)$, almost surely. Finally, using (\ref{Efirstp}) and the
definition of $B_\delta^{(n)}$, one can check inductively that
(\ref{Eonedominother}) holds. Namely, $\xi_n\leq\xi_n'$ is the
base of induction, and for $k\geq n$ either $-1\leq\xi_k\leq
\xi'_k$ (the last inequality is by induction hypothesis) in which
case $\Delta\xi'_k \geq\Delta\xi_k$ due to (\ref{Efirstp})
yielding $\xi_{k+1}\leq\xi'_{k+1} $, or $ \xi_k<-1$ and
$\xi'_k\geq\delta$ in which case on $B_\delta^{(n)}$ we have $
\xi_{k+1}<
\delta\leq\xi'_k=\xi'_{k+1}$, or finally $ \xi_k<-1$ and
$\xi_k\leq\xi'_k< \delta$ in which case again on $B_\delta
^{(n)}$ we have
$\xi_{k+1}=\xi_k+\Delta\xi_k \leq\delta\wedge(\xi'_k+\Delta
\xi_k)=\xi'_{k+1}$. Therefore,
\[
\mathbb{P}(\xi^* = \infty) \leq\mathbb{P}\bigl( \bigl(B_\delta^{(n)} \bigr)^c \bigr)+ \mathbb{P}\Bigl(\limsup_k {\xi'}^{(n,\delta)}_k =\infty\Bigr).
\]
We conclude that it suffices to show
%
\begin{equation}
\label{Esufficient}
\mathbb{P}\Bigl(\limsup_k {\xi'}^{(n,\delta)}_k =\infty\Bigr)=0
\end{equation}
for a fixed $\delta>0$ and each $n\geq1$.

Again by shift and scaling of space, and additional shift of time,
we can henceforth assume that $a=b=0$, and that (\ref{Efirstpeq})
holds. It is clear that if the process $\xi$ switches sign only
finitely many times then it either spends all but finitely many
units of time being nonnegative, in which case by the martingale
convergence theorem it converges, or it spends all but finitely
many units of time being nonpositive. On both events $\xi^*$ is
finite. It remains to prove the claim on the event $A^{\pm}$ where
$\xi$ switches sign infinitely often. In fact we will prove here a
stronger claim, namely that
%
\begin{equation}
\label{Estronger}
A^{\pm} \cap\{ \xi^*= 0\}{\ = A^{\pm} \cap\{ \xi^*\leq0\}} =
A^{\pm}
\qquad\mbox{almost surely.}
\end{equation}
The first identity above is clear from the definitions of
$A^{\pm}$ and $\xi^*$. Fix $\varepsilon>0$. For $n\geq1$, define the
process
\[
S_k^{(n)}:= \sum_{i=n}^{k-1} (\xi_{i+1}-\xi_i)1_{\{\xi_i\geq0\}},
\qquad k\geq n,
\]
with the convention $S_n^{(n)}=0$, and note that by assumption
(\ref{Efirstpeq}) on $\xi$, $S_\cdot^{(n)}$ is a martingale
started from $0$ at time $n$.

Due to Doob's maximal inequality we have
\[
\mathbb{P}\Bigl(\sup_{k\geq n} \big|S_k^{(n)}\big|>\varepsilon\Bigr)\leq\frac{4 \sum
_{k\geq n}
\mathbb{E}[(\xi_{k+1}-\xi_k)^2 1_{\{\xi_k\geq0\}}]}{\varepsilon^2}
\]
and in particular, due to (\ref{Ethirdp}), we can find $n_1\geq1$
such that this probability is smaller than $\varepsilon$, hence
%
\begin{equation}\label{EDoobi}
\mathbb{P}\Bigl(\sup_{k,j\geq n_1}\big|S_k^{(n_1)}-S_j^{(n_1)}\big|>2\varepsilon\Bigr)\leq2 \varepsilon.
\end{equation}
Consider $\xi$ on the event
\[
A^\pm\cap\Bigl\{\sup_{k,j\geq n_1} \big|S_k^{(n_1)}-S_j^{(n_1)}\big|\leq
2\varepsilon\Bigr\},
\]
and note that now the maximal value of $\xi$ on any excursion into
$[0,\infty)$ that begins after time $n_1$ cannot exceed
$\sup_{n\geq n_1} 1_{\{\xi_n<0 <\xi_{n+1}\}} \xi_{n+1}
+2\varepsilon\leq
o_{n_1}(1)+2\varepsilon$, where $o_{n_1}(1)\to0$, as $n_1\to\infty$.
Since $\varepsilon$ can be taken arbitrarily small, we obtain
(\ref{Estronger}).
\end{pf}

The above result (\ref{Estronger}) can be improved in the
following sense. Assume that $\xi$ satisfies
(\ref{Efirstp})--(\ref{Ethirdp}). Denote by $A_{a}^\pm$ the event
$\{\xi-a$ switches sign infinitely often$\}$.
\begin{lemma}
\label{Lswitchsign} On $A_{a}^\pm$, we have
\[
\xi^*\leq a+b,\qquad \mbox{a.s.}
\]
\end{lemma}
\begin{pf} We may assume again that $a=-1$ and $b=1$, and that
$\xi_0< -1$. Let $T_0=0$, and for $m\geq1$ let $T_m$ be the $m$th
downward crossing time of $-1$ by $\xi$. Note that on the event
$A_{-1}^\pm$, $T_m$ is finite\vspace*{-2pt} almost surely and that also $T_m\to
\infty$ as $m\to\infty$. It is clear how to generalize the
construction of ${\xi'}^{,(n,\delta)}$ from the proof of Lemma
\ref{Ljointexc} by replacing a fixed time $n$ by a stopping time
$T_m$, $m\geq0$. Of course, the construction extends only on the
event $\{T_m<\infty\}$, on the complement one can define the
process as identity $\delta$ (for example). We will henceforth
abbreviate ${\xi''}^{,(m,\delta)}\equiv{\xi'}^{,(T_m,\delta)}$.

Using (\ref{EDoobi}) and (\ref{Ethirdp}) one can easily check, as in
the proof of previous lemma,
that
\[
\lim_{m\to\infty} \sup_{k\geq T_m} {\xi_k''}^{,(m,\delta)}\leq
\delta.
\]
Indeed, the overshoots of ${\xi_k''}^{,(m,\delta)}$ across $\delta
$ are
becoming negligible as $m$ increases, and (\ref{Ethirdp}) controls
its fluctuations. In particular,
\[
\xi^* 1_{A_{-1}^\pm}\leq\Bigl(\lim_m\sup_{k\geq T_m}
{\xi''_k}^{,(m,\delta)}\Bigr) 1_{A_{-1}^\pm}\leq\delta.
\]
Since $\delta>0$ is arbitrary,
it follows that $\mathbb{P}(A_{-1}^\pm\cap\{\xi^* > 0\})=0$,
as claimed.
\end{pf}

\begin{rema}
\label{R:additconstr} We will sometimes consider a process $\xi$
adapted to the filtration $\mathcal{F}$, where the conditions
(\ref{Efirstp})--(\ref{Ethirdp}) apply up to additional
constraint. More precisely
\[
\mathbb{E}\bigl((\xi_{k+1}-\xi_k) 1_{\{\xi_k\geq a\}}|\mathcal{F}_k\bigr)
1_{E_k}\leq0,\qquad
\limsup_k
1_{\{\xi_k<a<\xi_{k+1}\}}
(\xi_{k+1} -a) 1_{E_k} \leq b,
\]
and
\[
\sum_{k} \mathbb{E}\bigl[(\Delta\xi_k)^2 1_{\{\xi_k\geq a\}\cap
E_k}\bigr]<\infty,
\]
where $E_k$ is an $\mathcal{F}_k$-measurable event. In such a
situation we
will (nonrigorously) state that $\xi$ satisfies (\ref
{Efirstp})--(\ref
{Ethirdp})
on $\bigcap_{k\geq n} E_k$ (for some large $n$) and conclude the
result of Lemma \ref{Ljointexc} on the same event.
The
corresponding rigorous formulation of this argument is to work
instead with the stopped process $\xi(T):=\{\xi_{k\wedge T},
k\geq n\}$, where a stopping time
\[
T:=\inf\{k\geq n\dvtx  1_{E_k}=0\}
\]
is defined precisely so that $\{T=\infty\} = \bigcap_{k\geq n} E_k$.
Then $\xi(T)$ satisfies the original (\ref{Efirstp})--(\ref{Ethirdp}),
and the asymptotics
of $\xi(T)$ and $\xi$ (as $k\to\infty$) match on
the event $\{T=\infty\}$.
\end{rema}

\section{Modified VRRW on a triangle}\label{sec_MT}

In this section we consider a modified VRRW (MVRRW) on a triangle.
Define $\tau_0^{(3)}=0$.
The transition probabilities of MVRRW are as for the VRRW on
the triangle, with one difference:
when the \textit{special vertex} $3$ is visited
for the $k$th time, at the stopping time
%
\begin{equation}
\label{Etauk} \tau_k^{(3)}\equiv   \tau_k:= \min\{t>
\tau_{k-1}\dvtx   X(t) = 3 \},\qquad k \geq1,
\end{equation}
its weight $Z(\tau_k,3)$ becomes $H(k)$ rather than $Z(\tau_k-1,3)+1$
[and for $t\in(\tau_k,
\tau_{k+1})$ we set $Z(t,3)=H(k)$].
Here we assume that the
sequence $H(k)$ is measurable with respect to $\mathcal{F}_{\tau_k}$, the
$\sigma$-algebra generated by the process up to time $\tau_k$,
that $H(1)\ge1$ and that for $k=0,1,2,\ldots$ the following property
holds:
\begin{eqnarray}\label{EH}
H(k+1)\ge H(k)+1.
\end{eqnarray}
Thus, the special vertex $3$ gets reinforced by a larger amount
than nonspecial vertices $1$ and $2$.

We study the above MVRRW with intention of applying it several
times in Section \ref{S:anaclike}. A typical application is in the
following context: suppose that the underlying graph is complete
graph on $d$ vertices where $d\geq4$. If one ``clumps together''
all but two of the vertices (say $i$ and $j$), then the VRRW
generates (with the appropriate time change) a MVRRW on a
triangle, where $i$ and $j$ correspond to $1$ and $2$, and the
clump corresponds to the special vertex $3$.

To simplify notation we will denote
\[
U(t):=Z(t,1), \qquad V(t) := Z(t,2) \quad\mbox{and}\quad W(t)=Z(t,3).
\]
The goal of this section is to show that $U(t) \asymp V(t)$. Before
stating the main result rigorously, we do some preliminary
comparisons and calculations.

First, observe that using elementary arguments (in particular,
{P\'{o}lya} urn-like transitions of the process, when viewed from the
special vertex $3$) one can show that for MVRRW both $U(t)\to\infty$
and $V(t)\to\infty$, almost surely. Similarly, it is easy to see
that it is impossible that after some finite time the particle
oscillates between nonspecial vertices $1$ and $2$. Hence
$W(t)\to\infty$, and $\tau_k<\infty$, for all $k$, almost surely.
Second, let us show that $W(t)$ cannot be too small with respect
to $U(t)+V(t)$ (which seems obvious but still requires a proof). Let
$\eta_n$, $n \geq0$ be the times of the successive visits to
vertices $1$ or $2$, that is
\[
\eta_{n+1}=\inf\bigl\{t>\eta_n\dvtx   X(t)\in\{1,2\} \bigr\}.
\]
Let $Y_n=W(\eta_n)$ and $X_n=U(\eta_n)+V(\eta_n)$. Then it is simple
to construct a coupling of
$(X_n,Y_n)$ with the urn
$(X_n',Y_n')$, featured in Theorem~\ref{thurn2} with $a=c=d=1$,
$b=0$, such that
\begin{eqnarray}\label{eqcouple}
X_n=X_n'\quad\mbox{and}\quad Y_n \geq Y_n' \qquad\mbox{for all } n.
\end{eqnarray}
This yields
\[
\liminf_{n\to\infty} \frac{Y_n}{X_n/\log X_n}\ge1.
\]
To simplify notation let
\[
\phi(x)=x/\log x.
\]
Then the above can be rewritten as
\begin{eqnarray*}
\liminf_{n\to\infty} \frac{W(\eta_n)}{\phi(U(\eta_n)+V(\eta
_n))}\ge1.
\end{eqnarray*}
Noting that in between the consecutive times $\eta_n$ the process
$W$ increases, while $U + V$ stays the same, we get
\begin{eqnarray}\label{3notsmall}
\liminf_{t\to\infty} \frac{W(t)}{\phi(U(t)+V(t))}\ge1.
\end{eqnarray}
Similarly, considering the process $(U(t),V(t),W(t))$ at times
when the MVRRW $X(t)$ visits vertex $1$ and comparing the
increments at vertices $1$ and $2$ [the former always increases by
$1$ while the latter increases by at least $1$ with probability at
least $V(t)/(U(t)+V(t))]$ we obtain that
\begin{eqnarray}\label{1notsmall}
\liminf_{t\to\infty} \frac{V(t)}{\phi(U(t))}\ge1,
\end{eqnarray}
and in a
symmetric way the symmetric result
\begin{eqnarray}\label{2notsmall}
\liminf_{t\to\infty} \frac{U(t)}{\phi(V(t))}\ge1.
\end{eqnarray}
To simplify notations further, recall (\ref{Etauk}), (\ref{EH})
and denote
\begin{eqnarray*}
U(\tau_k) = u,\qquad
V(\tau_k)= v,\qquad
W(\tau_k)=a=H(k),\qquad
n(k)=  n=u+v.
\end{eqnarray*}
We omit the index ``$k$'' from the notation in the forthcoming
argument, whenever not in risk of confusion. Relations
(\ref{3notsmall})--(\ref{2notsmall}) imply
(in a straightforward way) that for sufficiently
large $k$ we have
\begin{eqnarray}\label{Ebounduv}
u>\phi(v)/2, \qquad
v>\phi(u)/2\quad
&\Longrightarrow&\quad
\min\{u,v\}>\phi(n)/4 \quad
\mbox{and}\quad\nonumber
\\[-8pt]\\[-8pt]
&&{}\quad a>\phi(n)/2.\nonumber
\end{eqnarray}
At time $\tau_k+1$ the walk has to visit either site $1$ or~$2$, and moreover $\mathbb{P}(X(\tau_k+1)=1)=u/(u+v)$,
$\mathbb{P}(X(\tau_k+1)=2)=v/(u+v)$.

For $m\ge1$, consider the events
%
\begin{eqnarray}\label{EAmk}
\hspace*{10pt}A_m(k)&=&\bigl\{X(\tau_k+1)={\mathsf1},X(\tau_k+2)={\mathsf2},\nonumber
\\
&&\hspace*{4pt}{}X(\tau_k+3)={\mathsf1},X(\tau_k+4)={\mathsf2},\ldots,\\
&&\hspace*{4pt}{} X\bigl(\tau_k+(2m-1)\bigr)={\mathsf1},\mbox{ but }X(\tau_k+2m)={\mathsf3}\bigr\},\nonumber\\\label{EBmk}
\hspace*{10pt}B_m(k)&=&\{X(\tau_k+1)={\mathsf1},X(\tau_k+2)={\mathsf2},\ldots,X(\tau_k+2m-1)={\mathsf1},\nonumber\\[-8pt]\\[-8pt]
&&\hspace*{68pt}{}X(\tau_k+2m)={\mathsf2},\mbox{ but }X(\tau_k+2m+1)={\mathsf3}\}.\nonumber
\end{eqnarray}
Symmetrically define events $\bar{A}_m(k)$, $\bar{B}_m(k)$ where the
walker
starts the excursion away from vertex $3$ at vertex $2$, and on $\bar{A}_m(k)$
[resp.,~$\bar{B}_m(k)$]
it visits $2$ (resp.,~$1$) immediately before returning to $3$.
Note that $A_m,B_m$, $m \ge1$ are disjoint. On $A_m \cup B_m$,
during this excursion, vertex $1$ is visited exactly $m$
times, while vertex $2$ is visited $m-1$ times on $A_m$ and $m$
times on $B_m$. Symmetric statements apply to $\bar{A}_m$ and
$\bar{B}_m$. It is easy to see that
\begin{eqnarray*}
\mathbb{P}\biggl(\bigcup_m (A_m\cup B_m)\big| \mathcal{F}_{\tau_k}\biggr)
&=&\mathbb{P}\bigl(X(\tau_k+1)=1,\tau_{k+1}<\infty| \mathcal{F}_{\tau_k} \bigr)
\\
&=& \mathbb{P}\bigl(X(\tau_k+1)=1 |\mathcal{F}_{\tau_k}\bigr)\qquad \mbox{a.s.},
\end{eqnarray*}
since $\tau_{k+1}<\infty$, almost surely. Next observe that for
$m \ge1 $ (where an empty product is equal to $1$)
\begin{eqnarray*}
\mathbb{P}(A_m|\mathcal{F}_{\tau_k})=
\frac{u}{u+v} \prod_{j=0}^{m-2} \biggl( \frac{v+j}{v+j+a} \cdot
\frac{u+j+1}{u+j+1+a} \biggr) \frac{a}{a+v+m-1}
\end{eqnarray*}
and
\begin{eqnarray*}
\mathbb{P}(B_m|\mathcal{F}_{\tau_k})&=&
\frac{u}{u+v} \prod_{j=0}^{m-2} \biggl( \frac{v+j}{v+j+a} \cdot
\frac{u+j+1}{u+j+1+a} \biggr)
\\
&&\hspace*{44pt}{}\times\frac{v+m-1}{a+v+m-1}\frac{a}{a+u+m}.
\end{eqnarray*}
Now define
\[
C_m(k)\equiv C_m = \bigcup_{i=m}^\infty(A_i \cup B_i)
\]
to be the event that vertex $1$ is visited at least $m$ times
during the excursion (recall that there is dependence of $u,v,a$,
and hence of $A_m, B_m$, and $C_m$ on $k$). Then
\[
\mathbb{P}(C_m|\mathcal{F}_{\tau_k})= \frac{u}{u+v} \prod
_{j=0}^{m-2} \biggl(
\frac{v+j}{v+j+a} \cdot\frac{u+j+1}{u+j+1+a} \biggr).
\]
If we denote
\[
\l_u=\frac a{a+u},\qquad \l_v=\frac a{a+v}\quad \mbox{and}\quad
\nu=(1-\l_u)(1-\l_v)
\]
then, provided $m^2/u\ll1$ and $m^2/v\ll1$,
\begin{eqnarray}\label{EestiC}
&&\mathbb{P}(C_m{(k)}|\mathcal{F}_{\tau_k})\nonumber
\\
&&\qquad=\frac u{u+v}\cdot\nu^{m-1}\nonumber
\\[-8pt]\\[-8pt]
&&\qquad\quad{}\times\frac{ (1+\sfrac 0v ) (1+\sfrac1v )\ldots(1+\vfrac{m-2}v )}
{ (1+\afrac0{a+v} ) (1+\afrac1{a+v} )\ldots(1+\vafrac{m-2}{a+v} )}\nonumber
\\
&&\qquad\quad{}\times \frac{ (1+\sfrac1u ) (1+\sfrac 2u )\ldots(1+\vfrac{m-1}u )}{ (1+\afrac1{a+u} ) (1+\afrac 2{a+u} )\ldots(1+\vafrac{m-1}{a+u} )}\nonumber
\\
\label{EestiCa}
&&\qquad=\frac u{u+v}\cdot\nu^{m-1}\bigl(1+O(m^2/u)+O(m^2/v)\bigr).
\end{eqnarray}
Set $m=m(k)=\log^3 n(k)+1$,
then by (\ref{Ebounduv}) we have $m^2/u,\ m^2/v< 4\log^7
(n)/n=o(1)$. Similarly, by (\ref{Ebounduv}), we have
%
\begin{equation}
\nu=\frac{1}{(a/u+1)(a/v+1)}\le\frac{1}{(a/n+1)^2}\le\frac1{(1+
\afrac{1}{2\log n})^2},
\label{ECestprior}
\end{equation}
and so a straightforward calculus manipulation yields
\[
\nu^{m-1} \le n^{1- \log n}.
\]
Consequently,
%
\begin{equation}\label{ECest}
\hspace*{30pt}\mathbb{P}\bigl(C_{m(k)}(k)|\mathcal{F}_{\tau_k}\bigr)
=\mathbb{P}\bigl(C_{(\log n)^3 +1}|\mathcal{F}_{\tau_k}\bigr)< \nu^{m-1}\bigl(1+o(1)\bigr) \le\frac{1+o(1)}{n^{\log n-1}}.
\end{equation}
Therefore, by the Borel--Cantelli lemma,
%
\begin{equation}\label{ECs}
\mbox{only finitely many of $C_{m(k)}(k)$ occur, a.s.}
\end{equation}
If $m\le m(k)=\log^3 n+1$, then we can simplify the conditional
probabilities of $A_m$ and $B_m$ as follows:
%
\begin{eqnarray}\label{Eprobsim1}
\mathbb{P}(A_m|\mathcal{F}_{\tau_k} )&=&\frac u{u+v} \l_v \nu^{m-1}[1+O(\log^7 n/n)],\\
\mathbb{P}(B_m|\mathcal{F}_{\tau_k})&=&\frac u{u+v} \l_u(1-\l_v)\nu^{m-1}[1+O(\log^7n/n)],\\
\mathbb{P}(\bar{A}_m|\mathcal{F}_{\tau_k} )&=&\frac v{u+v} \l_u\nu^{m-1}[1+O(\log^7n/n)],\\\label{Eprobsim4}
\mathbb{P}(\bar{B}_m|\mathcal{F}_{\tau_k})&=&\frac v{u+v} \l_v(1-\l_u)\nu^{m-1}[1+O(\log^7 n/n)].
\end{eqnarray}
Now let
\[
\xi(t):= \frac{U(t)}{U(t)+V(t)}.
\]

\begin{lemma}
\label{Lliminf} We have
\[
\mathbb{P}\Bigl(\liminf_{t \to\infty} \xi(t)>0\Bigr)=1,
\]
and by symmetry $\mathbb{P}(\limsup_{t \to\infty} \xi(t)<1)=1$.
\end{lemma}
\begin{pf} It suffices to restrict attention to times $\tau_k$
since by (\ref{ECs}) the values of $\xi$ during the interval
$(\tau_k,\tau_{k+1})$ differ (asymptotically) from $\xi(\tau_k)$ by
at most
order $\log^3(U(\tau_k)+V(\tau_k))/(U(\tau_k)+V(\tau_k))$.
Recall that we abbreviate $V(\tau_k)=v$, $U(\tau_k)=u$, $n=u+v$.
In particular, $n\geq k+O(1)$ for each $k\geq1$, almost surely,
since between any two visits to site $3$, either site $1$ or~$2$ is visited at least once.

Define (recall the example in Section
\ref{S:polyaexample})
\[
\Xi(t)=\log\bigl(U(t)+V(t)\bigr)-\log\bigl(V(t)-1\bigr).
\]
We will estimate the drift of $\Xi$ (in the case where $v<n/3$,
hence $v<u/2$) by comparing our MVRRW setting to that of the
2-color {P\'{o}lya} urn. In the latter case, with probability
$u/(u+v)$ the new value is
\[
\mbox{{P\'{o}lya}}_\uparrow= \log(n+1)-\log(v-1),
\]
and with probability $v/(u+v)$ the new value is
\[
\mbox{{P\'{o}lya}}_\downarrow=\log(n+1)-\log(v).
\]
Thus, the drift increment of $\Xi$ under the law of the {P\'{o}lya} urn
is negative, since
\begin{eqnarray}\label{Esimpledrift}
\frac u{u+v}\log\frac{n+1}{v-1}
+\frac v{u+v} \log\frac{n+1}v -\log\frac{n}{v-1}<0,
\end{eqnarray}
see also Section \ref{S:polyaexample}.

Our goal is to bound the drift of $\Xi$ under the
modified VRRW law by its counterpart under
the {P\'{o}lya} urn process. Intuitively, this makes sense since the
shuttles pull the ratio $U/(U+V)$ closer to $1/2$, which
corresponds to even more negative drift of $\Xi$. Note that
\begin{eqnarray*}
&&\mathbb{E}(\Xi(\tau_{k+1}) | \mathcal{F}_{\tau_k} )
\\
&&\qquad= \sum
_{m=1}^\infty\biggl(
\mathbb{P}(A_m|\mathcal{F}_{\tau_k})\log{\frac{n+2m-1}{v+m-2}}
+\mathbb{P}(B_m|\mathcal{F}_{\tau_k})\log{\frac{n+2m}{v+m-1}}
\\
&&\qquad\quad\hspace*{22pt}{}+ \mathbb{P}(\bar A_m|\mathcal{F}_{\tau_k})\log{\frac{n+2m-1}{v+m-1}}
+\mathbb{P}(\bar B_m|\mathcal{F}_{\tau_k})\log{\frac{n+2m}{v+m-1}}
\biggr)
\\
&&\qquad= \bigl( \mathbb{P}(B_1|\mathcal{F}_{\tau_k}) +\mathbb{P}(\bar B_1|
\mathcal
{F}_{\tau_k}) \bigr)
\log{\frac{n+2}{v}}
\\
&&\qquad\quad{}+ \sum_{m=1}^\infty\biggl(
\mathbb{P}(A_m|\mathcal{F}_{\tau_k})\log\frac{n+2m-1}{v+m-2}
+\mathbb{P}(\bar A_m|\mathcal{F}_{\tau_k})\log{\frac
{n+2m-1}{v+m-1}} \biggr)
\\
&&\qquad\quad{}+\sum_{m=2}^\infty\bigl(
\mathbb{P}(B_m|\mathcal{F}_{\tau_k})
+\mathbb{P}(\bar B_m|\mathcal{F}_{\tau_k}) \bigr)\log{\frac{n+2m}{v+m-1}}
\\
&&\qquad={\mathsf I}+{\mathsf{II}}+{\mathsf{III}}.
\end{eqnarray*}
Then
\begin{eqnarray}
\label{EII}
{\mathsf {II}}\le\sum_{m=1}^\infty\biggl( \log{\frac{n+1}{v-1}}
\mathbb{P}(A_m|\mathcal{F}_{\tau_k})+ \log{\frac{n+1}{v}}\mathbb
{P}(\bar A_m|
\mathcal{F}_{\tau_k})
\biggr)
\end{eqnarray}
and
\begin{eqnarray}
\label{EIII}
{\mathsf {III}}\le\sum_{m=2}^\infty\biggl(\log
\frac{n+1}{v-1} \mathbb{P}(B_m|\mathcal{F}_{\tau_k})+ \log{\frac
{n+1}{v}}\mathbb{P}(\bar
B_m|\mathcal{F}_{\tau_k}) \biggr),
\end{eqnarray}
since for $m\ge2$ and $v<n/3$
\[
\frac{n+2m}{v+m-1}-\frac{n+1}{v}<0.
\]
Finally, since for $u>v$,
\begin{eqnarray*}
\mathbb{P}(B_1|\mathcal{F}_{\tau_k})&=&\frac{u}{n}\frac
{v}{v+a}\frac
{a}{a+u+1}\  >
\frac{v}{n}\frac{u}{u+a}\frac{a}{a+v+1}=\mathbb{P}(\bar
B_1|\mathcal{F}_{\tau_k}),
\end{eqnarray*}
we have
\begin{eqnarray}
\label{EcalcIst}
\hspace*{30pt}{\mathsf I}&=& \bigl(\mathbb{P}(B_1|\mathcal{F}_{\tau_k})+\mathbb{P}(\bar B_1|\mathcal{F}_{\tau_k})\bigr) \log\frac{n+2}{v}
\\
&=& \bigl(\mathbb{P}(B_1|\mathcal{F}_{\tau_k})-\mathbb{P}(\bar B_1| \mathcal{F}_{\tau_k})\bigr)
\log\frac{n+2}{v} + \mathbb{P}(\bar B_1|\mathcal{F}_{\tau_k}) 2\log{\frac{n+2}{v}}
\nonumber\\
&\le& \bigl(\mathbb{P}(B_1|\mathcal{F}_{\tau_k})-\mathbb{P}(\bar B_1|
\mathcal{F}_{\tau_k})\bigr)\log{\frac{n+2}{v}}\nonumber
\\
&&{}+ \mathbb{P}(\bar B_1|\mathcal{F}_{\tau_k})
\biggl(\log{\frac{n+1}{v-1}} + \log{\frac{n+1}{v}} \biggr)
\nonumber\\
&\le& \bigl(\mathbb{P}(B_1|\mathcal{F}_{\tau_k})-\mathbb{P}(\bar B_1|\mathcal{F}_{\tau_k})\bigr)\log{\frac{n+1}{v-1}}\nonumber
\\
&&{}+ \mathbb{P}(\bar B_1|\mathcal{F}_{\tau_k}) \biggl(\log{\frac{n+1}{v-1}} + \log{\frac{n+1}{v}} \biggr)
\nonumber\\
&=& \mathbb{P}(B_1|\mathcal{F}_{\tau_k}) \log{\frac{n+1}{v-1}} +\mathbb{P}(\bar B_1|\mathcal{F}_{\tau_k}) \log{\frac{n+1}{v}}.
\label{EcalcIen}
\end{eqnarray}
For the first inequality (the third line in the display) above we
use the fact that
\[
\biggl(\frac{n+2}{v} \biggr)^2 \leq\frac{(n+1)^2}{v(v-1)}
\qquad\mbox{whenever }v<\frac{n}{3}.
\]
Therefore,
\begin{eqnarray*}
{\mathsf I} + {\mathsf {II}} +{\mathsf {III}}
&\leq& \log{\frac{n+1}{v-1}}
\sum_{m=1}^\infty\bigl(\mathbb{P}(A_m|\mathcal{F}_{\tau_k})+\mathbb{P}( B_m|\mathcal{F}_{\tau_k})\bigr)
\\
&&{}+ \log{\frac{n+1}{v}} \sum_{m=1}^\infty\bigl(\mathbb{P}\bigl(\bar
A_m|\mathcal{F}_{\tau_k})+\mathbb{P}(\bar B_m|\mathcal{F}_{\tau_k})\bigr),
\end{eqnarray*}
and by noting
\[
\sum_{m=1}^\infty\bigl(\mathbb{P}(A_m|\mathcal{F}_{\tau_k})+\mathbb{P}(
B_m|\mathcal{F}_{\tau_k})\bigr)=u/(u+v)
\]
and
\[
\sum_{m=1}^\infty\bigl(\mathbb{P}(\bar A_m|\mathcal{F}_{\tau
_k})+\mathbb{P}(\bar
B_m|\mathcal{F}_{\tau_k})\bigr)= v/(u+v),
\]
we arrive to the following bound: provided $v<n/3$ (that is,
$v<u/2$), the drift increment of the $\Xi$ process under the
modified VRRW law is smaller than the expression on the left-hand side of
(\ref{Esimpledrift}). In particular, $\Xi$ has supermartingale
increments whenever its value is larger than $\log{4}$. It is
simple to check that $\Xi$ satisfies properties
(\ref{Efirstp})--(\ref{Ethirdp}) with $a= \log{4}$ [note that this
$a$ is different from $a\equiv a(k)$ above] and $b=0$ (any $b\geq
0$ would suffice). Namely, we have just verified (\ref{Efirstp}),
while (\ref{Esecondp}) is true since the steps
$\Xi(\tau_{k+1})-\Xi(\tau_k)$ are asymptotically of order at most
$\log^{4}(n)/n$, due to the lower bound (\ref{Ebounduv}) on $v$
and estimate (\ref{ECs}). Similarly, (\ref{Ethirdp}) holds since
\[
\bigg|\log\biggl(\frac{u+v+2m}{v-1+m} \biggr)- \log\biggl(\frac{u+v}{v-1} \biggr) \bigg|=
O \biggl(\frac{m}{v} \wedge\frac{u}{v} \biggr)= O \biggl(\frac{m \log{n}}{n}\wedge\log{n} \biggr),
\]
where the upper bound $u/v= O(\log{n})$ will be useful for atypically
large $m$.
Due to (\ref{ECestprior}), the above estimate implies the following
bound:
\begin{eqnarray}\label{Evaribnd}
&&\mathbb{E}\bigl(\bigl(\Xi(\tau_{k+1})- \Xi(\tau_k)\bigr)^21_{\{\Xi(\tau_k)\geq\log{4}\}}|\mathcal{F}_{\tau_k} \bigr)\nonumber
\\
&&\qquad\leq c \biggl[ \frac{\log^{8}{n}}{n^2} +\log^2{n} \times\mathbb{P}(C_{\log^3{n}+1}|\mathcal{F}_{\tau_k})\biggr]\\
&&\qquad\leq c \biggl( \frac{\log^{8}{n}}{n^2} +e^{-c' \log^2{n}} \biggr),\nonumber
\end{eqnarray}
where $c\in(0,\infty)$
and $c'\in(0,1)$ do not depend on $k$. Recall that $n\geq k$,
for all $k$, so the sequence (\ref{Evaribnd}) of upper bounds is
summable in $k$. Now Lemma \ref{Ljointexc} yields that
$\limsup_t \Xi(t)$ is finite almost surely, and this is
equivalent to saying that $\liminf_t \xi(t)$ is strictly positive,
almost surely.
\end{pf}

\section{Analysis on complete-like graphs}\label{S:anaclike}

We will denote by $\mathcal{G}=\mathcal{G}_d$ a
complete-like
graph of interest. Our main goal in this section is to prove the
following result leading to Theorem~\ref{th1}.
\begin{prop}\label{Pcomparison}
The VRRW on $\mathcal{G}$ satisfies:
\textup{(i)}
\[
\liminf_t \frac{Z(t,i)}{Z(t,j)}> 0 \qquad\mbox{a.s.},
\]
for any two different interior sites $i,j$.

\textup{(ii)} If $\ell_1,\ldots,\ell_r$ are the leaves attached to an interior
site $g$, then
\begin{eqnarray}
&&\biggl\{ \liminf_t \min_{i\neq g} \frac{\sum_{j\notin
\{i,g\}}Z(t,j)}{\sum_{j\neq i}Z(t,j)} > \delta\biggr\} \subset
\biggl\{\limsup_t \frac{(\sum_{j=1}^r
Z(t,\ell_j))^{1+\delta}}{\sum_{i\neq g}Z(t,i)}=0 \biggr\}\hspace*{-19pt} \nonumber
\\[-8pt]\\[-8pt]
\eqntext{\mbox{a.s.},}
\end{eqnarray}
where the sums above [except for $\sum_{j=1}^r Z(t,\ell_j)$] are
taken over
the interior sites only.
\end{prop}

In the following subsections we prove the above proposition,
treating several different cases separately. Property (ii) above
will be used in the proof of Theorem~\ref{th1}. It gives a priori
bounds on the total empirical frequency of the leaves, that
simplify the large deviations estimates relative to the
corresponding argument in~\cite{V2001} (see Section \ref{secproofs} for
details).

\subsection{Graphs with leaves at a single vertex}

We start by considering the simplest noncomplete graph from
the class of graphs described in the\break \hyperref[Introduction]{Introduction}. Here there are
three ``interior'' sites $1$, $2$ and $3$, forming a triangle, and
there is an additional leaf $\ell_1^3=\ell\sim3$. As in the study
of MVRRW we will denote $U(t)=Z(t,1)$, $V(t)=Z(t,2)$,
$W(t)=Z(t,3)$ and, moreover,
\[
L(t)=Z(t,\ell).
\]
Clearly, the process $(U,V,W)$, observed only at times $(\sigma_k)_{k
\ge0}$, where $\sigma_0=t_0$ (assume without loss of generality that
$X_{t_0}\in\{1,2,3\}$) and
\[
\sigma_k:=\min\bigl\{j> \sigma_{k-1}\dvtx X_j\neq X_{\sigma_{k-1}}, X_j \in\{1,2,3
\}\bigr\},\qquad k\geq1,
\]
has the law of $(Z(t,1),Z(t,2),Z(t,3))$ generated by the motion of
a particle according to a MVVRW with a special vertex $3$.
Therefore, Lemma \ref{Lliminf} insures that $U(t)\asymp V(t)$, or
equivalently, that both
%
\begin{equation}
\label{EUandV} \limsup_{t\to\infty} \frac{U(t)}{V(t)}\quad \mbox{and}\quad
\limsup_{t \to\infty}
\frac{V(t)}{U(t)}
\end{equation}
are finite random variables, almost surely. As in (\ref{Etauk}),
denote by $\tau_k^{(g)}$ the time of the $k$th successive visit to
site $g$, where $g\in\{1,2,3\}$.
Easy comparison of
$(L(\tau_k^{(3)}),U(\tau_k^{(3)})+V(\tau_k^{(3)}))$ with the
{P\'{o}lya} urn ensures preliminary estimate
%
\begin{equation}
\label{Eprelim} \limsup_k
\frac{L(\tau_k^{(3)})}{U(\tau_k^{(3)})+V(\tau_k^{(3)})}<\infty, \qquad
\mbox{a.s.}
\end{equation}
As we will soon see, $L(\tau_k^{(3)})\ll
U(\tau_k^{(3)})+V(\tau_k^{(3)})$ as a lower (random) power. First
note that for any $t$
\[
W(t)\leq U(t+1)+V(t+1)+L(t+1)+W(t_0),
\]
so that (\ref{EUandV}) and (\ref{Eprelim}) imply
%
\begin{equation}
\label{EWlessU} \limsup_t \frac{W(t)}{U(t)}<\infty\qquad \mbox{almost surely},
\end{equation}
and in turn that
%
\begin{equation}
\label{EUVorder} \min\Biggl\{\liminf_t \frac{U(t)}{t},\liminf_t
\frac{V(t)}{t} \Biggr\}>0\qquad \mbox{almost surely}.
\end{equation}
Given (\ref{EWlessU}), it is now plausible that $W$ has the same
asymptotic order as $U$,
since its increase is ``helped'' by the existence of the leaf $\ell$.
Soft arguments based on comparison with a generalized urn yield
\begin{eqnarray}\label{EUcompW}
\limsup_t \frac{\phi(U(t))}{W(t)}<\infty,
\end{eqnarray}
but not more, and comparison with the VRRW on the pure triangle
does not seem to be useful either in proving the complement to
(\ref{EWlessU}). However, the drift increment comparison argument
of Lemma \ref{Lliminf} is robust enough. Namely, denote by
${\widetilde{W}}$
the process that starts as ${\widetilde{W}}(t_0)=W(t_0)$, and that increases
by amount $1$ at time $t+1$ if $X(t)\in\{1,2\}$ and $X(t+1)=3$
(i.e.,~whenever the site $3$ is visited from another interior
site), and that otherwise remains unchanged. Then
%
\begin{equation}
W(t) = {\widetilde{W}}(t)+ Z(t,\ell)-Z(t_0,\ell)={\widetilde
{W}}(t)+ L(t)-L(t_0)
\label{EtildaW}
\end{equation}
in particular, ${\widetilde{W}}(t) \leq W(t)$ for all $t$. Consider the
process
%
\begin{equation}
\label{EXimod} \Xi(k):=\log\bigl(U\bigl(\tau_k^{(2)}\bigr)+{\widetilde{W}}\bigl(\tau_k^{(2)}\bigr)\bigr)
-\log\bigl({\widetilde{W}}\bigl(\tau_k^{(2)}\bigr)-1\bigr),\qquad k\geq1,
\end{equation}
adapted to the $\sigma$-field $\mathcal{F}_{\tau_k}$, $k\geq1$ where
$\tau_k\equiv \tau_k^{(2)}$. Let $u=U(\tau_k)$, $v=W(\tau_k)$,
${\tilde{v}}={\widetilde{W}}(\tau_k)$, $a=V(\tau_k)$, $n=u+{\tilde
{v}}$, and note that the
drift of $\Xi$ at time $k$ (provided $v<u/2$) is still less or
equal to expression (\ref{Esimpledrift}); in particular it is
negative, as we reason next. It is necessary to interchange
sites $2$ and $3$ in the definitions (\ref{EAmk}) and (\ref{EBmk})
and the rest of this argument. While the conditional probabilities
of $A_m,\bar{A}_m$, $m\geq1$ and $B_m,\bar{B}_m$, $m\geq2$ are different
in the current setting where $\ell$ exists, the estimates in
(\ref{EII}) and (\ref{EIII}) only concern the number of shuttles
$m$ between the two sites. Therefore,
\begin{eqnarray}\label{Ewithtildev}
&&\mathbb{E}\bigl(\Xi(k+1) | \mathcal{F}_{\tau_k} \bigr)\nonumber
\\
&&\qquad\leq \bigl(
\mathbb{P}(B_1|\mathcal{F}_{\tau_k}) +\mathbb{P}(\bar B_1|
\mathcal{F}_{\tau_k}) \bigr)
\log{\frac{n+2}{{\tilde{v}}}}\nonumber \\[-8pt]\\[-8pt]
&&\qquad\quad{}+
\sum_{m=1}^\infty\biggl( \log{\frac{n+1}{{\tilde{v}}-1}}
\mathbb{P}(A_m|\mathcal{F}_{\tau_k})+ \log{\frac{n+1}{{\tilde
{v}}}}\mathbb{P}
(\bar A_m|\mathcal{F}_{\tau_k})
\biggr)\nonumber\\
&&\qquad\quad{}+
\sum_{m=2}^\infty\biggl(\log
\frac{n+1}{{\tilde{v}}-1} \mathbb{P}(B_m|\mathcal{F}_{\tau_k})+
\log
{\frac{n+1}{{\tilde{v}}}}\mathbb{P}(\bar
B_m|\mathcal{F}_{\tau_k}) \biggr).\nonumber
\end{eqnarray}
Next observe that $\mathbb{P}(B_1|\mathcal{F}_{\tau_k})$ does not change
under the new law, since possible shuttles between site $3$ and
its leaf $\ell$ before the step from $3$ to another interior site,
do not influence the conditional law of this step. Finally,
observe that $\mathbb{P}(\bar{B}_1|\mathcal{F}_{\tau_k})$ is
smaller than
$(v/n)(u/(u+a))(a/(a+v+1))$ under the new law, since possible
shuttles between site $3$ and its leaf $\ell$ that happen before the
step from $3$ to $1$, make the probability of the move from $1$ to
$2$ smaller than $a/(a+v+1)$. Thus the estimates
(\ref{EcalcIst}) and (\ref{EcalcIen}) can be carried out verbatim.
Due to (\ref{Ewithtildev}), and the fact ${\tilde{v}}\leq v$, we obtain
\begin{eqnarray*}
\mathbb{E}\bigl(\Xi(k+1)| \mathcal{F}_{\tau_k} \bigr)
&\leq& \log{\frac{n+1}{{\tilde{v}}-1}}
\cdot\frac{u}{u+v}+ \log{\frac{n+1}{{\tilde{v}}}}\cdot\frac
{v}{u+v}\\
&\leq& \log{\frac{n+1}{{\tilde{v}}-1}} \cdot\frac{u}{u+{\tilde{v}}}+
\log{\frac{n+1}{{\tilde{v}}}}\cdot\frac{{\tilde{v}}}{u+{\tilde{v}}},
\end{eqnarray*}
as claimed. In order to apply Lemma \ref{Ljointexc}, it remains to
estimate the quantities in (\ref{Esecondp}) and (\ref{Ethirdp}).
Before doing so, we show that $L$ is a smaller power of $U+V$, and
therefore of $W$. So fix $\beta\geq1$ and consider again the times
$\tau_k^{(3)}$, $k\geq1$ of successive visits to site $3$. Note
that $\tau_k^{(3)}$ is different from $\sigma_k$ above, and from
$\tau_k\equiv \tau_k^{(2)}$ linked to the definition of $\Xi$.
Abbreviate
\[
L_k:= L\bigl(\tau_k^{(3)}\bigr),\qquad
U_k:= U\bigl(\tau_k^{(3)}\bigr),\qquad
V_k:= V\bigl(\tau_k^{(3)}\bigr),\qquad
W_k:= W\bigl(\tau_k^{(3)}\bigr)=k.
\]

Then, if $\delta\in(0,1)$, on
\[
P_k^\delta:= \biggl\{\frac{U_k}{U_k+W_k}\wedge\frac{V_k}{V_k+W_k}>\delta\biggr\},
\]
we have
\begin{eqnarray}\label{Eratiodrift}
\nonumber
\mathbb{E}\biggl(\frac{L_{k+1}^\beta}{U_{k+1}+V_{k+1}} \Big|\mathcal{F}_{\tau_k^{(3)}}\biggr)
&\leq&\frac{(L_k+1)^\beta}{U_k + V_k} \cdot\frac{L_k}{U_k + V_k+L_k} \\
&&{}+\frac{(L_k)^\beta}{U_k + V_k+1} \cdot\frac{(1-\delta
)(U_k+V_k)}{U_k +
V_k+L_k} \\
&&{}+  \frac{(L_k)^\beta}{U_k + V_k+2} \cdot
\frac{\delta(U_k+V_k)}{U_k + V_k+L_k}.\nonumber
\end{eqnarray}
Namely, either the walk visits the leaf $\ell$ at time
$\tau_k^{(3)}+1$ and steps back to site $3$ at time
$\tau_k^{(3)}+2=\tau_{k+1}^{(3)}$, or it visits $\{1,2\}$ at time
$\tau_k^{(3)}+1$, and given this, it revisits the same set at time
$\tau_k^{(3)}+2$ with probability larger than $\delta$.

Using (\ref{EUandV}) and (\ref{EWlessU}) one easily sees that
\begin{eqnarray}
\label{Efullprob} \mathbb{P}\Bigl(\lim_{\delta\searrow0} \liminf_k
P_k^\delta\Bigr)=1.
\end{eqnarray}
From now on we take $\delta$ small and think about the behavior
of the process $(L_k)^\beta/(U_k+V_k)$ on $\bigcap_{k\geq n_0}
P_k^\delta$,
where $n_0$ is a large finite integer.
\begin{rema}
The part (a) of the next lemma will not be used in the sequel of the current
argument; however its argument will be needed in the next section.
\end{rema}
\begin{lemma}
\label{LLwithU} \textup{(a)} Estimate (\ref{Eprelim}) and $\liminf_t
(U(t)\wedge V(t))/\phi(t) >0$ are already sufficient for
\begin{eqnarray}
\label{ELbUVone} \lim_t \frac{L(t)}{U(t) +V(t)}=0\qquad \mbox{a.s.}
\end{eqnarray}

\textup{(b)} On $\bigcap_{k\geq n_0} P_k^\delta$, for any $\beta<1 + \delta$
we have
that
\begin{eqnarray}\label{ELbUV} \lim_t \frac{L(t)^\beta}{U(t)
+V(t)}=0\qquad
\mbox{a.s.}
\end{eqnarray}
\end{lemma}
\begin{pf} (a) We need a slightly more precise estimate than
(\ref{Eratiodrift}). Namely, keeping track of which interior
vertex ($1$ or $2$) the walk visits first, one obtains that\looseness=1
\begin{eqnarray}\nonumber
\mathbb{E}\biggl(
\frac{L_{k+1}^\beta}{U_{k+1}+V_{k+1}} \Big|\mathcal{F}_{\tau_k^{(3)}}
\biggr)&\leq&
\frac{(L_k+1)^\beta}{U_k + V_k} \cdot\frac{L_k}{U_k + V_k+L_k} \\
\nonumber
&&{}+\frac{(L_k)^\beta}{U_k + V_k+1} \cdot\frac{U_k W_k}{(U_k +
V_k+L_k)(V_k+W_k)} \\
&&{}+\frac{(L_k)^\beta}{U_k + V_k+2} \cdot\frac{U_k V_k}{(U_k +
V_k+L_k)(V_k+W_k)} \\
\nonumber
&&{}+\frac{(L_k)^\beta}{U_k + V_k+1} \cdot\frac{V_k W_k}{(U_k +
V_k+L_k)(U_k+W_k)} \\
&&{}+ \label{Eratiodriftprec} \frac{(L_k)^\beta}{U_k + V_k+2} \cdot
\frac{V_k U_k}{(U_k + V_k+L_k)(U_k+W_k)}.\nonumber
\end{eqnarray}
The right-hand side in (\ref{Eratiodriftprec}) equals\vspace*{-1pt}
\begin{eqnarray}\label{Erhs}
\frac{L_k^\beta}{(U_k + V_k)} (1 + R_k ),
\end{eqnarray}
with
$\beta=1$, and with\vspace*{-1pt}
\begin{eqnarray*}
R_k&=& 1\big/(U_k+V_k+L_k)\\
&&\hspace*{8.7pt}{}\times\biggl\{1 -\biggl(\frac{U_k W_k}{(U_k+V_k+1)(W_k+V_k)}+ \frac{V_k
W_k}{(U_k+V_k+1)(W_k+U_k)}
\biggr) \\
&&\hspace*{27pt}{}-2
\biggl(
\frac{U_k V_k}{(U_k+V_k+2)(W_k+V_k)}+
\frac{U_k V_k}{(U_k+V_k+2)(W_k+U_k)}
\biggr) \biggr\}.
\end{eqnarray*}
The last expression equals to
\begin{eqnarray*}
- \biggl(
\frac{U_k V_k}{(U_k+V_k+2)(W_k+V_k)}+
\frac{U_k V_k}{(U_k+V_k+2)(W_k+U_k)}
\biggr) + O \biggl(\frac{1}{U_k+V_k} \biggr).
\end{eqnarray*}
Now due to hypotheses of part (a)
we conclude that $U_k+V_k \asymp k$ and
$U_k \wedge V_k \geq c k/\log{k}$ for some positive random $c$.
Hence the leading term above has absolute value larger than a term of
order $1/\log{k}$.
In particular, the process $L_k/(U_k+V_k)$ is a positive
super-martingale, so it converges
almost surely to a finite limit.
However, the limit must be $0$, since on the event $\lim_k
L_k/(U_k+V_k) >0$
the
drift increment above is of the order at least $1/(k\log{k})$, so the
drift would not be summable
otherwise.
In this way one can also see that the asymptotic order of $L_k$ may not
be of
the form $k/a_k$, if $a_k$ converge to infinity sufficiently slowly so
that $\sum_k
1/(k\log{k} \times a_k)=\infty$. The last observation will not be used
in the sequel.

(b) Note that on $\bigcap_{k\geq n_0} P_k^\delta$, for any $\beta<1 +
\delta$ we
have the same expression~(\ref{Erhs}) for the right-hand side in
(\ref{Eratiodriftprec}), except that now $R_k$ is smaller than
\[
\frac{1}{U_k+V_k+L_k} \biggl(\beta-(1-\delta)-2\delta+
O \biggl(\frac{1}{L_k} \biggr)+O \biggl(\frac{1}{U_k+V_k} \biggr) \biggr).
\]
This can be seen already from (\ref{Eratiodrift}), since
$(L_k+1)^\beta/L_k^\beta= 1+ \beta/L_k+ O(L_k^2)$. Consequently,
$R_k$ is
again negative for all sufficiently large $k$, and therefore
$L_k^\beta/(U_k+V_k)$ converges to a finite random quantity. In
particular, for any $\beta'<\beta$ the limit in (\ref{ELbUV}) is
$0$ on
the event $\bigcap_{k\geq n_0} P_k^\delta$, and due to (\ref{Efullprob}),
after letting $\delta\to0$, one obtains (\ref{ELbUV}), hence part
(ii) of Proposition \ref{Pcomparison} for the triangle ornamented
with a single leaf.
\end{pf}

In order to prove (\ref{Esecondp}) and (\ref{Ethirdp}) for the
process $\Xi$ from (\ref{EXimod}), we will
derive analogues to (\ref{ECest}) and (\ref{ECs}). The reader can
check that in the special case where the leaves are attached to $3$ only
(that is, no leaves are attached at $1$ or $2$), one does not need
(\ref{ELbUV}) to obtain sufficiently good estimates. Nevertheless,
we will soon consider the general case, hence doing the
calculations while accounting for (\ref{ELbUV}) will prove useful.

Due to Lemma \ref{LLwithU}(b) and (\ref{EUcompW}) and (\ref{EtildaW}), we
have $\{\bigcap_{k\geq n_0} P_k^\delta\}\subset\{{\widetilde
{W}}(t)\asymp
W(t)\}$, and therefore
\begin{eqnarray}\label{EUcomptW}
\biggl\{\bigcap_{k\geq n_0} P_k^\delta\biggr\} \subset\biggl\{\limsup_t \frac{\phi
(U(t))}{{\widetilde{W}}(t)}<\infty\biggr\}\qquad \mbox{almost surely.}
\end{eqnarray}
Suppose that $\beta>1$ and that $(p_k^m)_{m\geq1,k\geq1}$ is a
table of numbers in $(0,1)$ such that
\[
1-p_k^m\leq\frac{c(m,k)}{k^{1-1/\beta}},\qquad m,k \geq1,
\]
where, for each finite integer $s$,
%
\begin{equation}
\label{Egeomlikeass}
\limsup_k \max_{m\leq s} c(m,k) < \infty.
\end{equation}
Let $(G_k, k\geq0)$ be a
random process (adapted to a filtration $({\mathcal{H}}_k, k\geq0)$) taking
values in the
nonnegative integers, and assume that it satisfies
conditional ``geometric-like'' relations
%
\begin{equation}
\label{Egeomlike}
\mathbb{P}(G_k>m+1 | G_k > m,{\mathcal{H}}_{k-1})= 1 - p_k^{m+1},\qquad
m\geq0.
\end{equation}
Then
$\mathbb{P}(G_k>s| {\mathcal{H}}_{k-1})= \prod_{m\leq s} (1-p_k^m)
\leq
(\max_{m\leq s}
c(m,k))^s /k^{s(1-1/\beta)}$,
and therefore, under the assumption (\ref{Egeomlikeass}), we have
\begin{eqnarray}\label{Egeombd}
\lim_{j\to\infty}\mathbb{P}\biggl(\bigcap_{k\geq j} \{G_k\leq2 /(1-1/\beta
)\}\biggr)= 1.
\end{eqnarray}
Consider the behavior of VRRW on $\bigcap_{k\geq n_0} P_k^\delta$
and fix some $\beta\in(1,1+\delta)$. Following each time
$\tau_k^{(3)}=\sigma_{k'}$ when VRRW visits site $3$ from another
interior site, the particle will make a nonnegative (possibly
$0$) number $\widetilde{N}_{k}$ of shuttles to $\ell$ before
visiting the next (different) interior site at time
$\sigma_{k'+1}$. Note that $\widetilde{N}_{k}$ in fact stands for
$W(\sigma_{k'+1})-W(\sigma_{k'})=W(\sigma_{k'+1})-k$. Let $j$ be a
large integer. Since $W(\tau_k^{(3)})=W_k=k$ we have on
$\bigcap_{k\geq n_0} P_k^\delta$ that $U_k+V_k \geq2\delta
k/(1-\delta)$, and
due to (\ref{ELbUV}) that $L_k\leq k^{1/\beta}$, for all $k\geq j$
(with an overwhelming probability as $j\to\infty$). As a
consequence, one can construct a process $G$ satisfying
(\ref{Egeomlikeass}) and (\ref{Egeomlike}) [where $c(m,k)$ can be
taken as $2\delta/(1-\delta)$ for all $k\geq j$ and $m\leq s$, so the
$\limsup$ in (\ref{Egeomlikeass}) is bounded by $2\delta/(1-\delta)$]
such that $\widetilde{N}_{k}\leq G_k$ (note that $G$ is defined
for all $k$, but the coupling of $\widetilde{N}_k$ and $G_k$ is
necessary only for $k$ such that $\tau_k^{(3)}=\sigma_{k'}$). Due
to (\ref{Egeombd}), we conclude that
%
\begin{equation}
\label{EwtildeN}
\{\widetilde{N}_{k} \leq2 /(1-1/\beta)\}\qquad \mbox{for all sufficiently
large $k$},
\end{equation}
with an overwhelming probability on $\bigcap_{k\geq n_0} P_k^\delta$.

Therefore, one can redo the calculation (\ref{EestiC}), this time
writing instead of the third term an analogous
\begin{eqnarray}
\label{Ethirdnew}
\frac{ (1+\sfrac0v ) (1+\sfrac{s_1}v )\cdots(1+\sfrac{s_{m-2}}v )}{ (1+\afrac0{a+v} ) (1+\afrac{s_1}{a+v} )\cdots(1+\afrac{s_{m-2}}{a+v} )},
\end{eqnarray}
where $s_{i+1}-s_i\geq1$ and $s_{i+1}-s_i\leq2/(1-1/\beta)$ for all
$i$, and for all large $k$. The estimate (\ref{EestiCa}) holds as
before, with different
constants comprised in $O(m^2/u)+O(m^2/v)$. Together with
(\ref{EUVorder}), this immediately implies (\ref{ECest}) and
(\ref{ECs}), and thus (\ref{Esecondp}) and (\ref{Ethirdp}) for
$\Xi$, as at the end of the proof of Lemma \ref{Lliminf}. Note
that in this step we also make use of the preliminary estimate
(\ref{EUcomptW}).

The above reasoning applied on the event $\bigcap_{k\geq n_0} P_k^\delta$
only (see also Remark \ref{R:additconstr}), but due to
(\ref{Efullprob}) we conclude the following lemma.
\begin{lemma}
\label{LXimod}
\[
\limsup_t \Xi(t)<\infty\qquad \mbox{a.s.}
\]
\end{lemma}

As a consequence, $\liminf{\widetilde{W}}(t)/(U(t)+{\widetilde
{W}}(t))>0$, almost surely,
and since $W(t)\geq{\widetilde{W}}(t)$,
\[
\liminf W(t)/\bigl(U(t)+W(t)\bigr)>0\qquad \mbox{a.s.},
\]
completing the proof of Proposition~\ref{Pcomparison}(i) in the
special case of the graph with three interior vertices and one
leaf.

As the reader will quickly check, the proof above carries over to
any $\mathcal{G}$ with the same interior sites $\{1,2,3\}$ and finitely
many leaves $\{\ell_1,\ldots,\ell_r\}$, all attached to the interior
site $3$. Namely, for the purposes of the calculation in Lemmas
\ref{LLwithU} and \ref{LXimod} all the leaves can be combined into
one ``super-leaf'' so that, in particular, Proposition~\ref{Pcomparison} holds via the same argument.

Moreover, suppose that $\mathcal{G}$ has interior sites
$\{1,2,\ldots,d\}$, $d\geq4$, and finitely many leaves
$\{\ell_1,\ldots,\ell_r\}$, all attached to the interior site $d$. Let
the initial position $X(t_0)$ take value in $\{1,\ldots,d\}$,
almost surely. Fix two different sites $i,j\in\{1,\ldots,d-1\}$,
and define three classes
\begin{eqnarray}
\label{Eclasses}
C_1:=\{i\},\qquad C_2:=\{j\}\quad \mbox{and}\quad C_3:=\{1,\ldots,d\}\setminus\{i,j\}
\end{eqnarray}
of interior vertices.
Consider $S(t)= \sum_{h=1}^3 h 1_{\{X(t) \in C_{h}\}}$,
and a sequence of stopping times $\sigma_0:=t_0$,
\[
\sigma_k:=\min\{s > \sigma_{k-1}\dvtx  S(s)\neq S(\sigma_{k-1}),S(s)\neq 0\},\qquad
k\geq1.
\]
Note that the process
\begin{eqnarray}\label{EXprim}
X'\equiv\bigl(X'(k),k\geq0\bigr) = \bigl(S(\sigma_k), k \geq0\bigr)
\end{eqnarray}
is identical in law to the position process $X$ of a MVRRW, with a
special vertex $3$.
Indeed, $\{S(t)= h\}=\{X(t)\in C_h\}$, for $h=1,2,3$, and $(\sigma
_k)_{k\geq0}$
are the successive times when $X$ jumps from one class of interior
vertices to another.
Therefore, setting
\[
Z'(k,h):=\sum_{v\in C_h} Z(\sigma_k,v),\qquad h=1,2,3,
\]
it is simple to check that the transitions of $X'$ are driven by (\ref
{Etrans}),
with $X'$ (resp.,~$Z'$) replacing $X$ (resp.,~$Z$).
Moreover,
$Z'(k+1, 1)- Z'(k,1)$ [resp.,~$Z'(k+1, 2)- Z'(k,2)$] equals $1$ if $X'(k)=1$
(resp.,~$=$2), while $Z'(k+1,3)- Z'(k,3)=H(k)\geq1$ if $X'(k)=3$.
A careful reader will note that the measurability
requirement on $H$ (see the beginning of Section \ref{sec_MT})
necessitates considering $X'$ with respect to
stopped filtration $(\mathcal{F}_{\sigma_k})_{k\geq0}$ generated by $X$.
As before,
these observations ensure that $Z'(k,1)\asymp Z'(k,2)$ as $k\to\infty$.
Since $Z(t,i)=Z'(k,1)$ and $Z(t,j)=Z'(k,2)$, where $t\in[\sigma_k,
\sigma
_{k+1})$, we conclude
that $Z(t,i)$ and $Z(t,j)$ are asymptotically
comparable,
for all $i,j\in\{1,\ldots, d-1\}$, almost surely.
It is again easy to verify that
\[
\limsup_t \frac{Z(t,d)}{\sum_{i=1}^{d-1} Z(t,i)}<\infty\quad \mbox{and}\quad \limsup_t \frac{\phi(\sum_{i=1}^{d-1} Z(t,i))}{Z(t,d)}
<\infty,
\]
almost surely. Since the walk necessarily returns to $d$ after
each visit to a leaf, we have $L(t)\leq Z(t_d)+L(t_0)$, and
therefore by the first estimate above we conclude
\[
t=Z(t,d)+\sum_{i=1}^{d-1} Z(t,i) + L(t) = O\Biggl(\sum_{i=1}^{d-1} Z(t,i)\Biggr)\qquad
\mbox{almost surely}.
\]
This implies readily that $\sum_{i=1}^{d-1} Z(t,i) \asymp t$, and
therefore that $Z(t,1) \asymp t$ (or equivalently, $Z(t,i) \asymp
t$, $\forall i=1,\ldots,d-1$), almost surely. Again combine all
the leaves into a single super-leaf $\ell\sim d$. The
calculation of Lemma \ref{LLwithU}(b), for the process observed at
successive times $\tau_k^{(d)}$ of visit to site $d$, yields as
before Proposition \ref{Pcomparison}(ii). Finally, let
$U(t)=Z(t,1)$, $V(t)=\sum_{g=2}^{d-1} Z(t,g)$ and $W(t)=Z(t,d)$,
and consider the process at the successive times
%
\begin{equation}
\label{Etiprime}
\hspace*{20pt}\sigma_k':=\min\bigl\{j> \sigma_{k-1}'\dvtx X_j\neq X_{\sigma_{k-1}'}, X_j
\in\{2,\ldots,d-1\}\bigr\},\qquad k\geq1,
\end{equation}
of visit to the subset $\{2,\ldots,d-1\}$.
Set ${\widetilde{W}}(t_0)={\widetilde{W}}(t_0)$ and let
\[
{\widetilde{W}}(t):=W(t)-\bigl(Z(t,\ell)-Z(t_0,\ell)\bigr),\qquad t\geq t_0.
\]
Then the process $\Xi$ defined as in (\ref{EXimod}) (with $\sigma
_k'$ in
place of $\tau_k^{(2)}$) again
satisfies (\ref{Efirstp})--(\ref{Ethirdp}) with $a=\log{4}$ and
$b=0$, so Lemma \ref{LXimod} follows, implying Proposition~\ref{Pcomparison}(i) as before.

\subsection{General complete-like graphs with $d\geq3$}\label{S:gencomp}

Assume that we are given a general complete-like graph $\mathcal
{G}=\mathcal{G}_d$ from
\hyperref[Introduction]{Introduction}. Here the argument is somewhat more delicate, due to
the fact that we cannot anymore use the MVRRW to easily obtain
$Z(t,i)\asymp t$ for most (all but one) sites, which was essential
in applying Lemma~\ref{LLwithU}.

We start again by making some soft observations. If $\ell\sim
g$, then $Z(t,\ell)\leq Z(t+1,g)+Z(t_0,\ell)$ implies
that $t=\sum_{v\in V(\mathcal{G})} Z(t,v) \leq\sum_{i=1}^d (r_i +1)Z(t+1,i)
+O(1)$, and in
particular that
\begin{eqnarray}\label{Eearlyli} \liminf_t \sum_{i=1}^d
Z(t,i)/t>0,
\end{eqnarray}
almost surely.
Moreover, {P\'{o}lya}'s urn comparisons, as in Section~\ref{sec_MT},
imply that
\[
\sup_t Z(t,v)=\infty,\qquad v\in V(\mathcal{G}),
\]
and, for each $i$,
%
\begin{equation}
\limsup_t \frac{{\sum_{j=1}^{r_i}} Z(t,\ell_j^i)}{\sum_{g=1,g\neq
i}^d Z(t,g)}< \infty\qquad \mbox{almost surely.}
\label{Erationadditional}
\end{equation}
Here we recall that $\ell_j^i$, $j=1,\ldots,r_i,$ are the leaves attached
at the interior site $i$.
Soon we will see that the limit in
(\ref{Erationadditional}) is $0$. Since
%
\begin{equation}
\label{Einequaladdit}
Z(t,i)\leq{\sum_{j=1}^{r_i}} Z(t+1,\ell_j^i) +
\sum_{g=1,g\neq i}^d Z(t+1,g) + Z(t_0,i),
\end{equation}
after adding $\sum_{g=1,g\neq i}^d Z(t,g)$
to both sides, (\ref{Eearlyli}) and (\ref{Erationadditional}) yield
\begin{eqnarray}
\label{EZjointorder}
\hspace*{28pt}\liminf_t \sum_{g=1,g\neq i}^d Z(t,g)/t>0
\qquad\mbox{for each interior site } i, \mbox{ almost surely.}\hspace*{-12pt}
\end{eqnarray}
Without loss of generality assume that $X(t_0)\in\{1,\ldots,d\}$.
Moreover, as already noted, each visit to a leaf of $i$ is
immediately followed by a visit to $i$. Therefore, if $Z(0,i)>
\sum_{j=1}^{r_i} Z(0,\ell_j^{i})$, then
\begin{eqnarray}
\label{Esillyin} Z(t,i)> \sum_{j=1}^{r_i} Z(t,\ell_j^{i}),\qquad t\geq t_0,
\end{eqnarray}
and provided (\ref{Esillyin}) holds at some time $t$, it will
continue to hold at all later times. We claim that, for each
$i=1,\ldots,d$, (\ref{Esillyin}) holds starting from some finite
time. Indeed, due to~(\ref{EZjointorder}) the walk will almost
surely (eventually) make at least $(\sum_{j=1}^{r_i}
Z(0,\ell_j^i)-Z(0,i))^{+} + 1$ steps from $i$ to another interior
vertex, and this ensures (\ref{Esillyin}) upon the next return to
$i$. Starting from the finite (stopping) time at which
(\ref{Esillyin}) holds for all $i\in\{1,\ldots,d\}$, one can
compare (as in Section~\ref{sec_MT}) the process $(\sum_{g=1,g\neq
i}^d Z(\sigma_k,g),Z(\sigma_k,i))$, where $\sigma_k$ is the time of $k$th
return to the subset of sites $\{1,\ldots,d\}\setminus\{i\}$,
with the generalized urn $(X_k',Y_k')$ of Theorem~\ref{thurn2}
(again here $a=c=d=1$, $b=0$), so that $Z(\sigma_k,i)\geq Y_k'$ and
$\sum_{g=1,g\neq i}^d Z(\sigma_k,g)\leq X_k'$. In particular, for
each $i=1,\ldots, d$,
\begin{eqnarray}\label{Eapriest}
\liminf_t \frac{Z(t,i)}{\phi(\sum_{g=1,g\neq i}^d Z(t,g))}> 0
\quad\mbox{hence}\quad \liminf_t \frac{Z(t,i)}{\phi(t)}> 0\nonumber
\\[-8pt]\\[-8pt]
\eqntext{\mbox{almost surely.}}
\end{eqnarray}
Due to (the argument of) Lemma \ref{LLwithU}(a), estimates (\ref
{Erationadditional})
[namely, its consequence (\ref{EZjointorder})]
and (\ref{Eapriest}) are sufficient to conclude that almost surely,
for each $i=1,\ldots, d$,
\begin{eqnarray}
\label{Eleaf_inter}
\lim_t \frac{\sum_{j=1}^{r_i} Z(t,\ell_j^i)}{\sum_{g=1,g\neq i}^d Z(t,g)}=
\lim_t \frac{\sum_{j=1}^{r_i} Z(t,\ell_j^i)}{t}=0.
\end{eqnarray}
Indeed, the reader can quickly check that $\sum_{j=1}^{r_i} Z(t,\ell_j^i)$
[resp.,~$\sum_{g=1,g\neq i}^d Z(t,g)$], observed at the times of return
to $i$,
corresponds
to $L(t)$ [resp.,~$U(t)+V(t)$], observed at the times of return to $3$.
The possible presence of leaves at sites $g\neq i$, corroborates  inequality
(\ref{Eratiodrift}).

However, we wish to strengthen (\ref{Eleaf_inter}) to an analogue
of Lemma \ref{LLwithU}(b). In order to be able to recycle its
argument, it suffices to show that for any $i\neq g$, $i,g\in
\{1,\ldots,d\}$ we have
\begin{eqnarray*}\liminf_t \frac{\sum_{l=1,l\notin
\{i,g\}}^d Z(t,l)}{t}>0,
\end{eqnarray*}
or equivalently, that the third most
frequently visited interior site has positive asymptotic
frequency. Let $(Z_{(1)}(t),\ldots, Z_{(d)}(t))$ be the vector of
order statistics for $Z(t,g),g=1,\ldots,d$, and set
\begin{eqnarray*}
S(t)=Z_{(d)}(t),\qquad
P(t)= Z_{(d-1)}(t)\quad \mbox{and}\quad
R(t) =\sum_{j=1}^{d-2} Z_{(j)}(t).
\end{eqnarray*}
Clearly $S(t)\asymp t$, and due
to (\ref{EZjointorder}) also $P(t)\asymp t$. Moreover, due to
(\ref{Eleaf_inter}) it must be
\begin{eqnarray}\label{EPasy}
\liminf_t \frac{P(t)}{t}\geq\frac{1}{2(d-1)}.
\end{eqnarray}
Indeed, (\ref{Eleaf_inter}) implies that $\limsup_t S(t)/t \leq
1/2$, and hence, the identity
$S(t)+P(t)+R(t)+\sum_{i=1}^d\sum_{j=1}^{r_i}Z(t,\ell_j^i)\equiv t$
and (\ref{Eleaf_inter}) together imply\break $\liminf_t
(P(t)+R(t))/t\geq1/2$, and (\ref{EPasy}) in turn.

It suffices to show that $R$ is asymptotically comparable to $S+P$.
Let $a(t)=\min\{i\dvtx  Z_{(d)}(t)=Z(t,i)\}$ and $b(t)=\min\{i \neq a(t)\dvtx
Z_{(d-1)}(t)=Z(t,i)\}$.
Consider the process $\tilde{\eta}(t):=(S(t)+P(t))/R(t)$ at successive
times of visit to
the set $\{a(t),b(t)\}$.
Without risk of confusion,
let us denote by $(\tilde{\eta}_k, k\geq0)$ the process $\tilde
{\eta}$ viewed
only on this restricted collection of times.
\begin{lemma}
\label{Letatil}
$\limsup_k \tilde{\eta}_k < \infty$, almost surely.
\end{lemma}
\begin{pf} Let $\tau$ be the time of the $k$th visit to the set
of vertices $\{a(\cdot),b(\cdot)\}$. For concreteness suppose that
the current position $X(\tau)=b(\tau)$, the calculation below is
similar if $X(\tau)=a(\tau)$. Let $s,p,r$ denote the values of
$S(\tau),P(\tau),R(\tau)$, respectively, and let $l$ denote the
corresponding ``total leaf weight'' at $b(\tau)$. Without loss of
generality we may assume that $r\geq4(d-1)\geq4$. Assume in
addition that $s+p\geq2r$, or equivalently, that $\tilde{\eta}_k=
(s+p)/r\geq2$. Then, on $\{Z_{(d-1)}(\tau)> Z_{(d-2)}(\tau)\}$,
$\tilde{\eta}_{k+1}$ will either take value $(s+p+1)/r$ with
probability $(s+l)/(s+l+r)$, or a value smaller than
$(s+p+1)/(r+1)$ (here we use the fact that $s+p\geq2r$ and $r\geq
4$) with probability $r/(s+l+r)$. A careful reader will note that
this includes transitions that change values of $a$ or $b$. On the
opposite event $\{Z_{(d-1)}(\tau)= Z_{(d-2)}(\tau)\}$ it could be
that the particle jumps from $b(t)$ to another site with the same
frequency thus increasing $s+p$ by $1$ without changing $r$.
However, if
\begin{eqnarray}\label{Eestir}
r\leq\frac{\afrac{1}{3(d-1)}}{1-\afrac{1}{3(d-1)}} (p+s) \quad\Longrightarrow\quad
r \leq\frac{1}{3(d-1)} \tau,
\end{eqnarray}
then due to (\ref{EPasy}) we have $Z_{(d-2)}(\tau) < r \ll p$,
whenever $\tau$ is sufficiently large. In particular,
$\{Z_{(d-1)}(\tau)= Z_{(d-2)}(\tau)\}$ happens at most finitely
often, almost surely. Hence, provided $\tilde{\eta}_k \geq
3(d-1)\geq2$, the drift increment of $\tilde{\eta}$ is bounded by
\begin{eqnarray*}
\frac{1}{r} \cdot\frac{s+l}{s+r+l} - \frac{1}{r+1} \frac{s+p-r}{s+r+l},
\end{eqnarray*}
and since $r\geq4(d-1)$, it will be negative for all sufficiently
large $\tau$ due to (\ref{Eleaf_inter})--(\ref{Eestir}). It is particularly easy to check the other two
hypotheses of Lemma \ref{Ljointexc}. Namely, the absolute value of
the increment $\tilde{\eta}_{k+1}- \tilde{\eta}_k$ is of the order
$1/r=1/\sum_{g, g\neq a(\tau),b(\tau)} Z(\tau,g)$, so clearly
diminishing at the time instances when $\tilde{\eta}_k$ traverses
the threshold $3(d-1)$. Furthermore, due to (\ref{Eapriest}), the
sum of square increments is finite, a.s. The conclusion is now due
to Lemma \ref{Ljointexc}.
\end{pf}

It is easy to see that Lemma \ref{Letatil} implies $\liminf_t R(t)/t
>0$, and that this is equivalent
to having
\begin{eqnarray}
\label{EzaPdelta}
\liminf_t \min_{i,j=1}^d \frac{\sum_{g=1, g\notin\{i,j\}}^dZ(t,g)}{t} >0 \qquad\mbox{almost surely}.
\end{eqnarray}
In analogy to the setting of the previous subsection, for each
$g=1,\ldots, d$, define
\[
P_k^{\delta,g} := \biggl\{ \min_{i=1}^d \frac{\sum_{j=1, j\notin\{i,g\}}^d
Z(\tau_k^{(i)},j)}{\sum_{j=1,j\neq i}^d Z(\tau_k^{(i)},j)} \geq\delta\biggr\},
\]
where, as usual, $\tau_k^{(i)}$ is the $k$th return time to $i$.
The argument of Lemma \ref{LLwithU}(b) gives
\begin{eqnarray}
\label{ELLwithUanal}
\bigcap_{k\geq n_0} P_k^{\delta,g} \subset
\biggl\{\limsup_t \frac{(\sum_{j=1}^r
Z(t,\ell_j^{(g)}))^{\beta}}{\sum_{i\neq g}Z(t,i)}=0 \biggr\}
\end{eqnarray}
for any $\beta< 1+\delta$, and
this in turn yields Proposition \ref{Pcomparison}(ii).
Due to (\ref{EzaPdelta}), we have, moreover,
\begin{eqnarray}
\label{Efullproba}
\mathbb{P}\Biggl(\lim_{\delta\to0} \liminf_k \bigcap_{i=1}^d P_k^{\delta,i}\Biggr)=1.
\end{eqnarray}

Finally, consider two different interior sites $i$ and $j$, the
classes (\ref{Eclasses}) and the process $X'$ from
(\ref{EXprim}). In analogy to (\ref{EXimod}) and (\ref{Etiprime}),
for $g\in\{i,j\}$, define
\[
{\widetilde{Z}}(t,g):= Z(t,g) - \sum_{j=1}^{r_g}\bigl(Z(t,\ell_j^g)-
Z(t_0,\ell_j^g)\bigr),\qquad t\geq t_0.
\]
Then ${\widetilde{Z}}(t,g)\leq Z(t,g)$, $t\geq t_0$, $g\in\{i,j\}$, and,
moreover,
\begin{eqnarray}\label{Eshuttbd}
\Biggl\{ \liminf_k \bigcap_{i=1}^d
P_k^{\delta,i}\Biggr\}\subset\{{\widetilde{Z}}(t,j)\asymp Z(t,j),
{\widetilde{Z}}(t,i)\asymp
Z(t,i)\}\qquad\nonumber
\\[-8pt]\\[-8pt]
\eqntext{\mbox{almost surely.}}
\end{eqnarray}

Let $\sigma_k$
be the time of $k$th visit to class $C_3$ from $i$ or from $j$ (in
particular, not accounting for the steps from $C_3$ to itself, and
the steps from the leaves into $C_3$). Now consider
%
\begin{equation}
\label{EXimoda}
\hspace*{20pt}{\widetilde{\Xi}}(k):=\log\bigl({\widetilde{Z}}(\sigma
_k,i)+{\widetilde{Z}}(\sigma_k,j)\bigr) -
\log\bigl({\widetilde{Z}}(\sigma_k,j)-1\bigr),\qquad k\geq1.
\end{equation}
Fix $\delta\in(0,1)$ and $\beta<1+\delta$. The asymptotics
(\ref{ELLwithUanal}) ensures [see the discussion comprising
(\ref{Egeomlikeass})--(\ref{Egeombd})] the existence of a finite
$n_1$ such that with an overwhelming probability there are at most
$2/(1-1/\beta)$ repeated shuttles from $i$ (resp.,~$j$) to its leaves
following any step into $i$ (resp.,~$j$) from another interior site
that occurs during the time interval $(\sigma_k,\sigma_{k+1})$, for all
$k\geq n_1$.

We will show that a Doob--Meyer modification of the process
${\widetilde{\Xi}}$
still satisfies the properties (\ref{Efirstp})--(\ref{Ethirdp}) so
that again
\begin{eqnarray}\label{Eclosee} \limsup_k {\widetilde{\Xi}}(k)
<\infty \qquad\mbox{a.s. on } \liminf_k \bigcap_{i=1}^d P_k^{\delta,i}.
\end{eqnarray}
This is equivalent to
\[
\liminf_t \frac{{\widetilde{Z}}(t,j)}{{\widetilde{Z}}(t,i)}>0
\qquad\mbox{a.s. on }\liminf_k
\bigcap_{i=1}^d P_k^{\delta,i}.
\]
Due to (\ref{Efullproba}) and (\ref{Eshuttbd}) we can conclude
Proposition \ref{Pcomparison}(i).

Denote $u(k)\equiv u =Z(\sigma_k,i)$, ${\tilde{u}}(k)\equiv{\tilde{u}}
={\widetilde{Z}}(\sigma_k,i)$, $v(k)\equiv v =Z(\sigma_k,j)$,
${\tilde{v}}(k)\equiv
{\tilde{v}}={\widetilde{Z}}(\sigma_k,j)$, $n(k)\equiv n={\tilde
{u}}+{\tilde{v}}$ and $a(k)\equiv
a=\sum_{g\in C_3} Z(\sigma_k,g)$. In fact,
(\ref{Esecondp}) and (\ref{Ethirdp}) hold for ${\widetilde{\Xi}}$ as
in the case
of the graph with leaves at a single vertex only, using
(\ref{ELLwithUanal}) instead of Proposition \ref{Pcomparison}(ii). For (\ref{Efirstp}), note first that (cf.~also the next
lemma)
\[
\mathbb{P}(\bar B_1|\mathcal{F}_{\sigma_k})\leq\frac{v}{u+v}\cdot
\frac
{u}{a+u}\cdot
\frac{a}{a+v+1}\qquad \mbox{almost surely},
\]
since possible shuttles to leaves $\ell_1^j,\ldots,\ell_{r_j}^j$ can
only decrease the probability of return to class $C_3$ when
stepping out of $i$ into an interior site.
\begin{lemma}
\label{LbdonB_1} We have
%
\begin{eqnarray} \label{EbdonB_1}
\mathbb{P}(B_1|\mathcal{F}_{\sigma_k})\in\biggl[\frac{u}{u+v}\cdot\frac{v}{a+v}\cdot\frac{a(1- \varepsilon(k))}{a+u+1},
\frac{u}{u+v}\cdot\frac{v}{a+v}\cdot\frac{a}{a+u+1} \biggr]\hspace*{-40pt}\nonumber
\\[-8pt]\\[-8pt]
\eqntext{\mbox{almost surely},}
\end{eqnarray}
where $\varepsilon(k)$ is $\mathcal{F}_{\sigma_k}$-measurable
nonnegative random
variable, such that
on\break $\bigcap_{k\geq n_0} P_k^{\delta,i}$,
\[
\varepsilon(k) = O \biggl(\frac{(a+v)^{1/\beta}}{a+u} \biggr)\qquad \mbox{almost surely}.
\]
\end{lemma}
\begin{pf}
Recall that on $B_1$ the particle steps from a site in the class $C_3$
to $i$,
next does a
certain number $N(k;u)$
(possibly $0$) of shuttles to the leaves $\ell_1^i,\ldots,\ell_{r_i}^i$
before a step
to $j$,
and finally, does a number
(possibly $0$) of shuttles to the leaves $\ell_1^j,\ldots,\ell_{r_j}^j$
before stepping back to $C_3$.
It is now simple to check that
\[
\varepsilon(k)= \frac{u+v}{u} \mathbb{E}\biggl[
1_{\{X(\sigma_k+1)=i\}}
\mathbb{E}\biggl( \frac{N(k;u)}{a+u+N(k;u) + 1}
\Big|\mathcal{F}_{\sigma_k}, X(\sigma_k+1)=i \biggr)
\Big|\mathcal{F}_{\sigma_k} \biggr],
\]
so it suffices to show (recall that $v<u/2$)
\[
\mathbb{E}\biggl( \frac{N(k;u)}{a+u+N(k;u)} \Big|\mathcal{F}_{\sigma_k},
X(\sigma_k+1)=i \biggr)
\leq C \frac{(a+v)^{1/\beta}}{a+u}\qquad \mbox{almost surely},
\]
for some finite constant $C$.
Let $q\equiv q(k):= \sum_{j=1}^{r_i} Z(\sigma_k,\ell_j^i)\equiv
\sum
_{j=1}^{r_i} Z(\sigma_k+1,\ell_j^i)$
be the total weight of the leaves attached to $i$ at time $\sigma_k$ (that
is, $\sigma_k+1$).
Our calculation is based on the same
reasoning as the discussion
comprising (\ref{Egeomlikeass})--(\ref{Egeombd}); however,
the expectation bound is simpler,
since the random variable $N(k;u)/(a+u+N(k;u))$ of interest is bounded
by 1.
Namely,
$\mathbb{P}(N(k;u) \geq2q|\mathcal{F}_{\sigma_k}, 1_{\{X(\sigma
_k+1)=i\}
}) \leq
\mathbb{P}(N(k;u)\geq q +1|\mathcal{F}_{\sigma_k}, 1_{\{X(\sigma
_k+1)=i\}})=
\frac{q}{a+u+q}$,
and therefore
\[
\mathbb{E}\biggl( \frac{N(k;u)}{a+u+N(k;u)} \Big|\mathcal{F}_{\sigma_k},
X(\sigma_k+1)=i \biggr)
\leq\frac{2q}{a+u+2q} +\frac{q}{a+u+q} \leq\frac{3q}{a+u}.
\]
The very last term is bounded by ${C (v+a)^{1/\beta}}/(a+u)$,
provided $q \leq C(v+a)^{1/\beta}$, which happens eventually on
$\bigcap_{k\geq n_0} P_k^{\delta,i}$, almost surely.
\end{pf}

Note that almost surely on $\{v<u/2\}$
%
\begin{equation}
\label{Ebdepsk}
\frac{(a+v)^{1/\beta}}{a+u} = O \biggl(\frac{1}{(a+u)^{1-1/\beta}} \biggr)=
O \biggl(\frac{1}{(\sigma_k)^{1-1/\beta}} \biggr),
\end{equation}
where we used (\ref{EZjointorder}) for the last estimate.
Due to the fact $\mathbb{P}(B_1|\mathcal{F}_{\sigma_k})
+\varepsilon(k)
\geq\mathbb{P}(\bar
B_1|\mathcal{F}_{\sigma_k})$ the calculations
(\ref{EcalcIst}) and (\ref{EcalcIen}) can be modified to yield
\begin{eqnarray*}
&&\bigl(\mathbb{P}(B_1|\mathcal{F}_{\sigma_k})+\mathbb{P}(\bar B_1|
\mathcal{F}_{\sigma_k})\bigr) \log\frac{n+2}{{\tilde{v}}}
\\
&&\qquad\leq \mathbb{P}(B_1|\mathcal{F}_{\sigma_k}) \log{\frac
{n+1}{{\tilde{v}}-1}} + \mathbb{P}(\bar B_1|\mathcal{F}_{\sigma_k}) \log{\frac{n+1}{{\tilde{v}}}}
\\
&&\qquad\quad{}+ \varepsilon(k) \biggl(\log\frac{n+1}{n+2} +\log\frac{{\tilde{v}}}{{\tilde{v}}-1} \biggr).
\end{eqnarray*}
Denote
\[
r(k):=
\varepsilon(k) \biggl(\log\frac{n+1}{n+2} +\log\frac{{\tilde
{v}}}{{\tilde{v}}-1} \biggr) 1_{\{v < u/2\}}.
\]
We therefore obtain
\begin{eqnarray}\label{Eforlater1}
&&\mathbb{E}\bigl({\widetilde{\Xi}}(k+1) -{\widetilde{\Xi}}(k)|
\mathcal{F}_{\tau_k} \bigr)\nonumber
\\
&&\qquad\leq \log{\frac{n+1}{{\tilde{v}}-1}} \cdot\frac{u}{u+v}+ \log
{\frac{n+1}{{\tilde{v}}
}}\cdot\frac{v}{u+v} -\log{\frac{n}{{\tilde{v}}-1}}+r(k)
\nonumber\\
&&\qquad\leq \frac{1}{u+v} \biggl[\frac{u+v}{{\tilde{u}}+{\tilde{v}}}-\frac
{v}{{\tilde{v}}} \biggr]+O \biggl( \frac
{1}{{\tilde{v}}\cdot n} \biggr)+r(k)\\
&&\qquad= \frac{1}{u+v}\cdot\frac{{\tilde{v}}u - {\tilde{u}}v}{({\tilde
{u}}+{\tilde{v}}){\tilde{v}}}+O \biggl( \frac
{1}{{\tilde{v}}\cdot n} \biggr)+r(k)\nonumber\\
&&\qquad\leq \frac{1}{u+v}\cdot\frac{u({\tilde{v}}-v)+v(u-{\tilde
{u}})}{({\tilde{u}}+{\tilde{v}}){\tilde{v}}}+O \biggl(
\frac{1}{{\tilde{v}}\cdot n} \biggr)+r(k) \label{Eforlater2}\\
&&\qquad=:  \tilde{r}(k)\nonumber,
\end{eqnarray}
where for the second inequality we develop (recall $n = {\tilde
{u}}+{\tilde{v}}$)
\begin{eqnarray*}
\log\biggl(\frac{{\tilde{u}}+{\tilde{v}}+1}{u+v+1} \biggr)
- \log\biggl(\frac{{\tilde{u}}+{\tilde{v}}}{u+v} \biggr)
\quad\mbox{and}\quad
\log\biggl(\frac{v}{{\tilde{v}}} \biggr) - \log\biggl(\frac{v-1}{{\tilde{v}}-1} \biggr)
\end{eqnarray*}
via Taylor's expansion up to quadratic order terms.
Lemma \ref{LbdonB_1}, jointly with (\ref{Eapriest}),
(\ref{ELLwithUanal}) and (\ref{Ebdepsk}), implies that, on
$\bigcap_{k\geq n_0} \bigcap_{i=1}^d P_k^{\delta,i}$, $D_\infty:=
\sum_{l=1}^\infty\tilde{r}(l)$ is a finite random variable,
almost surely. Now observe that on $\{D_\infty\leq
K\}=\bigcap_{k\geq1} \{\sum_{l=1}^k \tilde{r}(l)\leq K\}$, the
process
\[
{\widetilde{\Xi}}' := \biggl({\widetilde{\Xi}}(k) - \sum_{l\leq k-1}
\tilde{r}(l),\ k\geq0 \biggr)
\]
satisfies (\ref{Efirstp})--(\ref{Ethirdp}) with $a = \log{4} + K$
and $b=0$.
Indeed, as in the previous section, one can argue that (\ref{EwtildeN})
holds for both
shuttles to the leaves attached at $i$ and at $j$ on
$\bigcap_{k\geq n_0} \bigcap_{i=1}^d P_k^{\delta,i}$. Hence one can redo the
calculation (\ref{EestiC}), where this time the third term is replaced by
(\ref{Ethirdnew}), and the second one by an analogous expression.
Due to Lemma \ref{Ljointexc}, $\limsup_t
{\widetilde{\Xi}}'(t)<\infty$, thus $\limsup_t {\widetilde{\Xi
}}(t)\leq\limsup_t {\widetilde{\Xi}}'(t) +
K<\infty$ on $\{D_\infty\leq K\}$, almost surely. By taking $K$
arbitrarily large we obtain (\ref{Eclosee}).

\subsection{\texorpdfstring{Proof of Theorem \protect\ref{th1}}{Proof of Theorem 1}}\label{secproofs}

For a fixed $\varepsilon>0$ define events
\begin{eqnarray*}
A(t)=A_{\varepsilon}(t)= \biggl\{\min_{i=1,\ldots,d} \frac{Z(t,i)}t \ge
\varepsilon
\mbox{ and }\max_{i=1,\ldots,d} \frac{\sum_{j=1}^{r_i}Z(t,\ell
_j^i)}t\le
t^{-\varepsilon} \biggr\}.
\end{eqnarray*}
Let
\begin{eqnarray*}
C_{\varepsilon}= \Biggl\{\exists T\dvtx  \bigcap_{t=T}^{\infty}A_{\varepsilon
}(t) \mbox{ occurs} \Biggr\}.
\end{eqnarray*}
\begin{prop}
\label{Propone}
We have
$C_{\varepsilon} \subseteq\{\pi_{\infty}=\pi_{\mathsf{unif}}\}$,
almost surely.
\end{prop}

\begin{pf} The argument is effectively a copy of that for
Theorem~1 in~\cite{V2001}.
The only difference is that now
the event $C_{\varepsilon}$ guarantees that the events $E(k)$ defined on
page 73 of~\cite{V2001} occur for all large enough $k\ge K$ (see
\cite{V2001}, formula (3.1)).
Observe that $\varepsilon_*$ in the definition of $E_2'(k)$ might
need to be
chosen quite large, yet this does not cause difficulties in applying
the argument.
Indeed, $\varepsilon_*$ does not need to satisfy
\cite{V2001}, formulas (3.23) and (3.24), since we can skip step 5 of \cite{V2001}---in
the current setting it is already covered by our estimates in previous sections,
hence included in the event $C_{\varepsilon}$.
Consequently (see \cite{V2001}, pages~73--74, for the definition
of $\gamma(k)$ and $k_0$), we have that, whenever $k_0\ge K$,
\begin{eqnarray*}
\mathbb{P}( \pi_{\infty}=\pi_{\mathsf{unif}}| C_{\varepsilon})
&\ge& \mathbb{P}\bigl(\pi_{\infty}=\pi_{\mathsf{unif}}|
C_{\varepsilon
},E(k_0) \bigr) \mathbb{P}( E(k_0) | C_{\varepsilon})\\
&=&\mathbb{P}\bigl(\pi_{\infty}=\pi_{\mathsf{unif}}| C_{\varepsilon},E(k_0)
\bigr)\ge\prod_{k=k_0+1}^{\infty}
\bigl(1-\gamma(k)\bigr),
\end{eqnarray*}
which, since $\sum_k \gamma(k) < \infty$,
can be made arbitrarily close to $1$ by choosing sufficiently large
$k_0$.
\end{pf}

\begin{pf*}{Proof of Theorem~\ref{th1}} Let
\[
\xi_{ij}:=\liminf_{t\to\infty} \frac{Z(t,i)}{Z(t,j)}
\]
and $\tilde C_n= \{\min_{i,j: i\ne j} \xi_{ij}>\frac
1n \}$. Proposition~\ref{Pcomparison}(i) implies that
$\mathbb{P}(\bigcup_{n=1}^{\infty}
\tilde C_n )=1$, or equivalently,
%
\begin{eqnarray}\label{eqCn}
\lim_{n\to\infty} \mathbb{P}(\tilde C_n)=1.
\end{eqnarray}
On the other hand, by part (ii) of Proposition~\ref{Pcomparison} and some easy
algebra, we have $\tilde C_n \subset C_{\afrac1{nd}}$.
The claim now follows from Proposition~\ref{Propone}
and~(\ref{eqCn}).
\end{pf*}

\subsection{Case $d=2$}\label{S:casetwo}

In this section, we briefly discuss a somewhat singular case,
where the number of leaves attached to the two ``interior''
vertices $1$ and $2$ influences the qualitative asymptotic
behavior of the corresponding VRRW.

Namely, if $r_1=r_2=0$, we have trivially (deterministically) $\pi
_{\infty}\to
\pi_{\mathsf{unif}}$, in accordance with Theorem \ref{th1}. However,
if $r_1>0$
and $r_2=0$ then site $2$ becomes qualitatively equal to any leaf
of $1$, and easy (multi-color {P\'{o}lya} urn) arguments show that
$Z(t,1)/t\to1/2$, while $Z(t,2)/t\to\alpha/2$, where $\alpha$ is
a continuous random variable taking values in $[0,1]$. In
particular, here $\pi_{\infty}\not\to\pi_{\mathsf{unif}}$.
Finally, the most
interesting case is when $r_1 \cdot r_2 >0$. By combining as usual
all the leaves attached to the same interior vertex into a single
super-vertex, we can assume $r_1=r_2=1$. Then abbreviating
\[
U(t)=Z(t,1),\qquad V(t)=Z(t,2),\qquad L(t)= Z(t,\ell_1^1),\qquad R(t)= Z(t,\ell_1^2),
\]
one can easily check that $U(t) \asymp V(t) \asymp t$ as
$t\to\infty$. Moreover, the process $L/(L+V)$ is a supermartingale
when observed at times of successive visits to vertex $1$. The
symmetric statement holds for the process $R/(R+U)$. Due to the
nonnegative supermartingale convergence, the limits
\[
\xi_L:= \lim_{t\to\infty} \frac{L(t)}{L(t)+V(t)},\qquad
\xi_R:= \lim_{t\to\infty} \frac{R(t)}{R(t)+U(t)},
\]
both exists, almost surely. Comparison with the {P\'{o}lya} urn
implies $ \mathbb{P}(\xi_L=1)=\mathbb{P}(\xi_R=1)=0. $ Using
comparison with urns
featured in Theorem \ref{thurn1}, one realizes that $\{\xi_L
>0\}\subset\{\xi_R =0\}$, almost surely, and moreover that $R(t)
= o(t^{1/a})$ for any $a\in(1, 1/\xi_L)$. The same statement
holds with $L$ and $R$ interchanged. Clearly, $\pi_{\infty}\not\to
\pi_{\mathsf{unif}}$
on $\{\xi_L >0\}\cup\{\xi_R >0\}$.

The results of~\cite{V2001}, Theorem 1.1, indicate that each
$\{\xi_L
>0\}$ and $\{\xi_R >0\}$ happen with positive probability; however,
we do not have an argument for $\mathbb{P}(\{\xi_L
>0\}\cup\{\xi_R >0\})=1$.

Using the process ${\widetilde{\Xi}}$ from (\ref{EXimoda}), and the reasoning
analogous (but simpler to that)
of Section \ref{S:gencomp} we obtain for $\beta> 1$
%
\begin{equation}
\label{Emayhelp}
\{ L(t) = O(t^{1/\beta})\} \subset\{ \xi_R > 0\}.
\end{equation}

\section{Consequences for $d$-partite graphs with leaves}\label
{Sec:generalize}

Assume $d\geq3$. The following graph $\widetilde{\mathcal{G}}\equiv
\widetilde{\mathcal{G}}_d=({\widetilde{V}}_d,\widetilde{E}_d)$, featured
in~\cite{V2001} as an example of\vspace*{1pt} a trapping subgraph for VRRW. It
is a generalization of $\mathcal{G}_d$ from the \hyperref[Introduction]{Introduction}, where
${\widetilde{V}}$ is
partitioned into $d+1$ equivalence classes $V_1,V_2,\ldots,V_d,B$.
The classes $V_i$, $i=1,\ldots,d$ are called the \textit{generalized
vertices}, and satisfy
the following two ($d$-partite structure) properties:
\begin{longlist}[(ii)]
\item[(i)] if $x,y \in V_i$, for some $i\in\{1,\ldots,d\}$, then $x \not\sim y$;
\item[(ii)] if $x\in V_i$ and $y\in V_j$ for two different $i,j\in\{1,\ldots
,d\},$ then $x \sim y$.
\end{longlist}
Moreover,
$B=\bigcup_{i=1}^d B_i$, where $B_i$ contains the ``leaves'' of $V_i$,
$i\in\{1,\ldots,d\}$,
\begin{longlist}[(iii)]
\item[(iii)] if $x\in B$ then there exists a unique $i\in\{1,\ldots,d\}$
such that $x \sim y$ for at least one $y\in V_i$.
\end{longlist}
Let $X$ be a VRRW on $\widetilde{\mathcal{G}}_d$. Then $X'$ defined by
\begin{eqnarray*}
X'(t) &=& \cases{
i,&\quad $X(t)\in V_i, i=1,\ldots,d$,\cr
\ell_i,&\quad $X(t)\in B_i, i=1,\ldots,d,$
}
\\
Z'(t,i)&:=&\sum_{x \in V_i} Z(t,x),\qquad
Z'(t,\ell_i):=\sum_{y \in B_i} Z(t,y),\qquad t\geq t_0,
\end{eqnarray*}
is very closely related to VRRW on graph $\mathcal{G}_d$ with
$r_1=\cdots= r_d=1$.
Namely, the only difference is that
on $\{X'(t)=i\}$ (that is, on $\{X(t)\in V_i\}$)
some of the weight $Z'(t,\ell_i)$ may not be accounted for when
computing the probability of the step to $X'(t+1)$, since $X(t)$
may equal $x\in V_i$ that is not connected to all the leaves in $B_i$.

Our methodology of Sections \ref{sec_MT} and \ref{S:anaclike} carries
over to the current setting and
we obtain the almost sure convergence of
local time frequencies for $X'$
to $\pi_{\mathsf{unif}}$ defined for $\mathcal{G}_d$.
Moreover, as in Proposition \ref{Pcomparison}, the leaves $\ell
_1^1,\ldots
, \ell_d^1$
are asymptotically visited a lower power order of times compared to the
interior vertices.

This translates to the following almost sure behavior of the VRRW
on $\widetilde{\mathcal{G}}_d$: the asymptotic proportion of time
spent in
$V_i$ is $1/d$ for each $i\in\{1,\ldots,d\}$, while the number of
visits to $B$ up to time $t$ is of the order $t^\alpha$, for some
random $\alpha$ such that $\mathbb{P}(\alpha\in(0,1))=1$.

We end this discussion with the following observation.
If $x, y \in V_i$, for some $i\in\{1,\ldots,d\}$, then
%
\begin{equation}
\lim_{t\to\infty}\frac{Z(t,x)}{Z(t,y)}\in(0,1)
\qquad\mbox{almost surely.}
\label{EZxZy}
\end{equation}
Note that if $B_i=\varnothing$, (\ref{EZxZy}) is a trivial
consequence of the {P\'{o}lya} urn convergence (see Section
\ref{S:polyaexample}). Namely, in this case the returns to class
$V_i$ can happen only from $\bigcup_{j\neq i} V_j$ and they clearly
have the (multi-color) {P\'{o}lya} urn distribution. To see
(\ref{EZxZy}) if $B_i\neq\varnothing$, first note that as before
one can use simple coupling with the urn of Theorem~\ref{thurn2}
to obtain preliminary estimates
\begin{eqnarray}\label{Eeasas}
\liminf_{t\to\infty} \frac{Z(t,x)}{\phi(Z(t,y))}\ge1\qquad \forall
x,y \in V_i.
\end{eqnarray}
Let $L(t) =\sum_{i=1}^d Z'(t,\ell_i)$ count the
visits to all the leaves combined. Due to the observations made
two paragraphs above, we have that $\mathbb{P}(\bigcup_{\beta>1}
G_\beta
)=1$, where
$G_\beta:=\{Z'(t,i) \to1/d, L(t)= O(t^{1/\beta})\}$. The asymptotics
of $Z'(\cdot,i)$, combined with (\ref{Eeasas}), now imply that
\begin{equation}\label{Eeasas1}
\hspace*{20pt}\bigcap_{x\in V_i}\{Z(t,x) \geq\phi(t)/(2
|V_i|) \}\qquad \mbox{for all sufficiently large $t$, almost surely.}
\end{equation}

Assume WLOG that $X(t_0)\in\bigcup_{j\neq i} V_j$, let $\tau_0=t_0$
and for $k\geq1$ let $\sigma_k:=\inf\{t> \sigma_{k-1} \dvtx  X(t-1)\in V_i,
X(t)\in\bigcup_{j\neq i} V_j\}$ be the $k$th time of return to
$\bigcup_{j\neq i} V_j$ from the class $V_i$. Let
\begin{eqnarray*}
{\widetilde{Z}}(t,x)&:=& {\widetilde{Z}}(t-1,x) + 1_{\{X(t-1)\in\bigcup
_{j\neq i} V_j, X(t)=x\}},
\\
{\widetilde{Z}}(t,y)&:=& {\widetilde{Z}}(t-1,y) + 1_{\{X(t-1)\in\bigcup
_{j\neq i} V_j, X(t)=y\}},
\qquad t\geq t_0,
\end{eqnarray*}
counts the visits to $x$ and $y$, respectively, made from interior
points \textit{exclusively} (due to definition of $\widetilde
{\mathcal{G}}$,
these points are necessarily contained in generalized vertices
different from $V_i$). Note that $0\leq Z(t,x)-{\widetilde{Z}}(t,x)
\leq
L(t)$, so that
\begin{eqnarray}\label{EasympZZpr} \bigcap_{t\geq t_0}
\bigcap_{x\in V_i} \biggl\{\bigg |\frac{{\widetilde{Z}}(t,x)}{Z(t,x)}-
1 \bigg|\leq\frac{L(t)}{Z(t,x)} \biggr\} \qquad\mbox{almost surely.}
\end{eqnarray}
Due to (\ref{Eeasas1}), we conclude that $Z(t,x)/{\widetilde
{Z}}(t,x)\to
1$ on $G_\beta$, and by letting $\beta\searrow1$ that
$Z(t,x)/{\widetilde{Z}}(t,x)\to1$, almost surely. Therefore, in order
to show
(\ref{EZxZy}) it suffices to prove
\begin{eqnarray}\label{Ejustneed}
\liminf_{t\to\infty} \frac{{\widetilde{Z}}(t,x)}{\sum_{y\in V_i}
{\widetilde{Z}}(t,y)} =
\limsup_{t\to\infty} \frac{{\widetilde{Z}}(t,x)}{\sum_{y\in V_i}
{\widetilde{Z}}(t,y)} >
0\qquad \forall x \in V_i.
\end{eqnarray}
Define an ``analogue'' of
(\ref{EXimoda})
\[
{\widetilde{\Xi}}(k):=\log\biggl({\widetilde{Z}}(\sigma_k,x)+\sum_{y\in
V_i\setminus\{x\}} {\widetilde{Z}}(\sigma_k,y)
\biggr) -
\log\bigl({\widetilde{Z}}(\sigma_k,x)-1 \bigr),\qquad k\geq1,
\]
%
and note that estimates (\ref{Eeasas})--(\ref{EasympZZpr})
ensure that (on each $G_\beta$) ${\widetilde{\Xi}}$
is a supermartingale up to a summable drift. In particular, it is converging to a finite (random) limit. This setting is
quite similar to
that mentioned at
the very end of Section \ref{S:casetwo}, as the estimates are
simpler than those of (\ref{Eforlater1}) and (\ref{Eforlater2}) due
to the following fact: there is no extra term $r(k)$ in
(\ref{Eforlater1}) in the current setting, since there are no
direct ``shuttles'' from $x$ to $y$ on the interval
$(\sigma_k,\sigma_{k+1}]$, indirect ``communication'' of $x$ and $y$ via
a common leaf is atypical---its occurrence is accounted for by
the differences $Z(t,x)-{\widetilde{Z}}(t,x)$, $Z(t,y)-{\widetilde
{Z}}(t,y)$, that are
both bounded by $L(t)$. Letting $\beta\searrow1$ establishes
(\ref{Ejustneed}). Let $Z_m(t)$ count the number of visits to
site $m$ up to time $t$ for VRRW on five (or fewer, at least
three) points $\{-2,-1,0,1,2\}$. Then the process
$(Z(t,x),Z(t,y))$ can be closely matched (coupled) to the process
$(Z_{-1}(t),Z_1(t))$ on the event $\{Z_{-2}(t)= O(t^{1/\beta_1}),
Z_2(t)= O(t^{1/\beta_2})\}$, where $\beta_1,\beta_2$ are two random
quantities strictly greater than $1$. The ``middle point'' $0$
corresponds to $\bigcup_{j\neq i} (V_j \cup B_j)$, while the
``boundary'' $-2$ (resp.,~$2$) corresponds to the set of leaves in
$B_i$ connected to $x$ (resp.,~$y$). Recall once again the process
${\widetilde{\Xi}}$ from (\ref{EXimoda}) and note that we are in the
situation
of type (\ref{Emayhelp}) where ${\widetilde{\Xi}}$ will be a supermartingale
up to a summable drift, and, moreover, where
${\widetilde{Z}}_{-1}(t)/Z_{-1}(t)\to1$ and ${\widetilde
{Z}}_1(t)/Z_1(t)\to1$. This
implies that $ \lim_{t\to\infty}\frac{Z_{-1}(t)}{Z_1(t)}\in
(0,1), \mbox{ almost surely,} $ hence (\ref{EZxZy}).

\section{Speed of convergence}\label{S:speed}

We first show a preliminary statement, which can be
viewed as a refinement of Proposition~3.2, page~80 in~\cite{V2001}.
\begin{lemma}\label{lemmaclaim1}
Suppose that we are given a sequence $(\eta_k)_{k\geq1}$ such
that for some $\varepsilon>0$ we have
\begin{eqnarray}
\hspace*{20pt}0\le\eta_k\le1-\varepsilon\quad\mbox{and}\quad
\eta_{k+1}\le\eta_k \biggl[1-\frac{C(1-\eta_k)} k \biggr]+\frac
{D}{k^{1+\tilde\beta}}
\qquad\forall k\ge k_0, \label{eqstatclaim}\hspace*{-12pt}
\end{eqnarray}
where $C>0$, $D>0$, and $\tilde\beta\in[0,1]$. Then $\limsup_{k\to
\infty}
\eta_k h(k) <\infty$,
where
\begin{eqnarray*}
h(k)= \cases{
k^{\tilde\beta}, &\quad \mbox{if }$\tilde\beta<C$,\cr
k^{\tilde\beta}/\log k, &\quad\mbox{if }$\tilde\beta=C$,\cr
k^{C}, &\quad\mbox{if }$\tilde\beta>C$.
}
\end{eqnarray*}
\end{lemma}
\begin{pf} First of all, let us show that $\eta_k\to0$. Indeed,
fix a positive $\tilde\varepsilon<\min\{C\varepsilon, \tilde
\beta\}$, and
suppose that
\begin{eqnarray}\label{eqetabound}
\eta_k\le\frac A{k^{\tilde\varepsilon}}
\end{eqnarray}
for some $A>0$. Then
\begin{eqnarray*}
\eta_{k+1}&\le& \frac A{k^{\tilde\varepsilon}} \biggl(1-\frac
{C\varepsilon} k \biggr)+\frac
{D}{k^{1+\tilde\beta}}
\\
&=&\frac A{(k+1)^{\tilde\varepsilon}}-\frac{A(C\varepsilon-\tilde
\varepsilon
)-Dk^{\tilde\varepsilon-\tilde\beta}-\Theta(k^{-1})}{k^{1+\tilde
\varepsilon}}
\le\frac A{(k+1)^{\tilde\varepsilon}},
\end{eqnarray*}
provided $A$ and $k$ are sufficiently large. We
obtain by induction that (\ref{eqetabound}) holds for all large $k$.
Therefore, one can, in fact, assume that $\varepsilon$ in
(\ref{eqstatclaim}) is arbitrarily close to $1$.
Hence, if
$\tilde\beta<C$, we can set $\tilde\varepsilon=\tilde\beta$
and, assuming that $\varepsilon
\in(0,1)$ is
sufficiently large
so that $C\varepsilon> \tilde\varepsilon$, we obtain (\ref
{eqetabound}) for any $A$
larger than
$D/(C\varepsilon- \tilde\varepsilon)=D/(C\varepsilon- \tilde
\beta)$.
This implies the claim of the lemma in the case $\tilde\beta<C$.

From now on assume $\tilde\beta\ge C$.
The above arguments imply that
for $\tilde\varepsilon=2C/3$, we have $\eta_k\le Ak^{-\tilde
\varepsilon}$, for all
large $k$
and some $A<\infty$, hence
\[
\eta_{k+1}\le\eta_k \biggl[1-\frac{C} k \biggr]+\frac{C\eta_k^2} k +\frac
{D}{k^{1+\tilde\beta}}\le
\eta_k \biggl[1-\frac{C} k \biggr]
+\frac{\bar D}{k^{1+\bar\beta}},
\]
where $\bar\beta=\min\{\tilde\beta,4C/3\}$ and $\bar D=D+A^2C$.
If
\[
\mu_k=\eta_k{k^{C}},
\]
then the last estimate together with Taylor's expansion of $(k+1)^C$
about $k$ yields
\begin{eqnarray*}
\mu_{k+1}&\le&\frac{\mu_k (k+1)^C}{k^C} \biggl[ 1- \frac{C} k \biggr] +\frac{\bar D(1+\Theta(1/k))}{k^{1+\bar\beta-C}}
\\
&\le&\mu_k \biggl[1-\frac{C(1+C)}{2k^2}+\Theta(k^{-3}) \biggr]+\frac{2\bar D}{k^{1+\bar\beta-C}}.
\end{eqnarray*}
By summing over $k$, this immediately implies
$\limsup_k \mu_k<\infty$ if $\tilde\beta>C$ (that is, $1+\bar
\beta-C>1$)
and and $\limsup_k \mu_k /\log k<\infty$ if $\tilde\beta=C$,
finishing the proof of the lemma.
\end{pf}

\begin{pf*}{Proof of Theorem~\ref{th2}}
Denote by
\begin{eqnarray*}
\eta(t):=1- d \min_{j=1,\ldots,d}\frac{Z(t,j)}{t} \in[0,1]
\end{eqnarray*}
another measure of distance between the empirical occupation measure
$\pi(t)=(Z(t,1)/t,\ldots,Z(t,d)/t)$ and $\pi_{\mathsf{unif}}=(1/d,\ldots,1/d)$.
Due to Theorem~\ref{th1}, we have
$\sum_j \pi_j(t)=1-o(1)$, so $\eta(t)/d\le
\Vert\pi(t)-\pi_{\mathsf{unif}}\Vert(1+o(1))\le\eta(t)$. Thus it
suffices to study the asymptotic behavior of $\eta(t)$.

Fix some constants $m>1$ and $\beta\in(0,(m-1)/2)$, and let
$\nu=\frac{m-1}2-\beta>0$.
Now consider VRRW at times $t_k=k^m$, set $N_k=t_{k+1}-t_k$ and
$\alpha_j^{(k)}=Z(t_k,j)/t_k,$ $j\in\{1,\ldots, d\}$, $k\in\mathbb
{N}$ (here
we use notations similar to those in the proof of Theorem~1
in~\cite{V2001}; also in order to simplify expressions we will often
omit the superscript ${}^{(k)}$ on $\alpha$'s).
Define events
\[
D_t(\varepsilon):= \bigcap_{i=1}^d \biggl\{\frac{Z(t,i)}{t}\in\biggl(\frac
1{d}-\varepsilon,\frac
1{d}+\varepsilon\biggr) \biggr\},\qquad t\geq t_0,
\]
and note that Theorem \ref{th1} can be rephrased as
\begin{eqnarray}
\label{Econvergence}
\hspace*{20pt}P\biggl(\forall\varepsilon\in(0,1/d) \mbox{ there is } K=K(\varepsilon
)<\infty\mbox{ s.t.} \bigcap_{k\geq K} D_k(\varepsilon) \mbox{ occurs}\biggr)= 1.
\end{eqnarray}
Fix some small positive $\varepsilon<1/d$. Due to (\ref
{Econvergence}) we can assume
from now on that $\min_j \alpha_j^{(k)}\ge\varepsilon$.

It is simple to check that if we were to ``freeze'' the
configuration at time $t_k$, ignore the visits to the leaves and let
the VRRW evolve as a Markov chain on state space $\{1,\ldots, d\}$
with transition probabilities
specified by the weights $(\alpha_j^{(k)})_{j=1}^d$ [or equivalently,
by $(Z(t_k,j))_{j=1}^d$], then this Markov chain would have
its reversible measure proportional to $(\alpha_1^{(k)}(1-\alpha
_1^{(k)}),\ldots,
\alpha_d^{(k)}(1-\alpha_d^{(k)}))$. As in the proof of \cite{V2001},
Theorem 1,
one uses the large deviation estimates (\ref{LDPg}) and (\ref{Hpro1})
to see that the number $N_{k:i}$ of visits to vertex $i$ during
$[t_k,t_{k+1})$ concentrates about its ``almost'' expected
value (i.e., the expectation according to the above frozen measure)
\begin{eqnarray}\label{EfrozenLLN}
\frac{\alpha_i(1-\alpha_i)}{\sum_{j=1}^d \alpha_j(1-\alpha
_j)} \times N_k=
\frac{\alpha_i(1-\alpha_i)}{1-\sum_{j=1}^d \alpha_j^2} \times N_k.
\end{eqnarray}
More precisely, let
\begin{eqnarray}\label{eqEkdef}
\hspace*{20pt}E_k&=& \bigl\{\mbox{simultaneously for all $i\in\{1,\ldots,d\}$,
the quantity $N_{k:i}$}\nonumber
\\[-8pt]\\[-8pt]
&&\hspace*{2pt}{} \mbox{ does not differ from (\ref{EfrozenLLN}) by more
than $k^{\vfrac{m-1}2+\nu} \asymp k^{\nu} \sqrt{N_k} $}
\bigr\}. \nonumber
\end{eqnarray}
Then (see \cite{V2001}, display~(3.16), page~76),
\begin{eqnarray*}
P(E_k^c)\le\gamma_k':= {\mathsf {Const}}_1(d) \exp(- {\mathsf {Const}}_2
(\varepsilon,d)
k^{2\nu} ),
\end{eqnarray*}
so we have $\sum_k
\gamma_k'<\infty$.
Therefore only finitely many $E_k^c$ occur.
Consequently, a.s. there is a $k_0=k_0(\omega)$ such that $\bigcap
_{k\ge
k_0} E_k$ occurs.
From now on, we will implicitly
assume that $k\ge k_0$.

We next recall that VRRW may also visit the leaves between
times $t_k$ and $t_{k+1}$. We already know from
Proposition~\ref{Pcomparison} that $\max_i \sum_{j=1}^{r_i}
Z(t,\ell_j^i)\le t^{1-\varepsilon'}$ for some $\varepsilon'>0$.
Let us now
strengthen this statement.
\begin{lemma}\label{lemrefinenew}
Let $L(t,i):=\sum_{j=1}^{r_i} Z(t,\ell_j^i)$ be the total cumulative
weight of all the leaves attached to $i$ at time $t$, where
$i\in\{1,\ldots,d\}$. Then, if $r_i>0$, for any $\delta>0$ we have
\begin{eqnarray}\label{eqd-1dnew}
\mathbb{P}\biggl(\liminf_{t\to\infty}\frac{L(t,i)}{t^{\afrac1{d-1}-\delta}}=\infty\biggr)=1
\end{eqnarray}
and (trivially if $r_i=0$)
\begin{eqnarray}\label{eqd-1d}
\mathbb{P}\biggl(\limsup_{t\to\infty}\frac{L(t,i)}{t^{\afrac1{d-1}+\delta}}=0 \biggr)=1.
\end{eqnarray}
\end{lemma}
\begin{pf} We will prove only the first part of the statement,
since the second one follows by an analogous argument.

As usual, let $\tau_k^{(i)}$ be the $k$th return time to the interior
vertex $i$.
Define
$X_k':=\sum_{g\neq i} Z(\tau_k^{(i)},g)$ and $Y_k':= L(\tau_k^{(i)},i)$.
Due to Theorem \ref{th1} and some simple algebra, the statement of the
lemma is
equivalent to the following claim: for any $\delta>0$ we have
\[
\limsup_{k\to\infty}\frac{X_k'}{(Y_k')^{{d-1}+\delta}}=0\qquad
\mbox{almost surely.}
\]
Recall (\ref{Econvergence}). Without loss of generality we
observe the process $(X', Y'):=((X_k', Y_k'), k\geq k_1)$, where
$\tau_{k_1}^{(i)}\geq K$ for some large finite $K$. In the spirit
of Remark~\ref{R:additconstr}, we will modify the VRRW and in
this way the process $(X', Y')$ (note, however, that here the
construction is slightly more complicated since we cannot simply
``truncate'' the process upon exiting the event of ``good
behavior''). Fix some small $\varepsilon>0$, and define
\[
D_t'(\varepsilon):= \bigcap_{i=1}^d \biggl\{\frac{Z(t,i)}{\sum_{j=1}^d Z(t,j)}\in
\biggl(\frac1{d}-\varepsilon,\frac1{d}+\varepsilon\biggr) \biggr\},\qquad t\geq t_0.
\]
Due to (\ref{Econvergence}) and Proposition \ref{Pcomparison}(ii)
we have that
\begin{eqnarray}
\label{EPDprime}
P\biggl(\bigcap_{k\geq K} D_k'(\varepsilon)\biggr)\to1\qquad\mbox{as } K\to\infty.
\end{eqnarray}
Define
\[
T_\varepsilon(K)\equiv T_\varepsilon:=\inf\{l> K \dvtx
D_l'(\varepsilon)\mbox{ does not occur}\}.
\]
If $K>2/\varepsilon$, it is easy to see that $D_{l-1}'(\varepsilon
)\subset
D_{l}'(3\varepsilon/2)$ for $l\geq K$, so
\begin{eqnarray}\label{Eintbd}
\{T_\varepsilon<\infty\}\subset\bigcap_{K\leq l\leq T_\varepsilon}
D_l'(3\varepsilon/2)
\qquad\mbox{almost surely.}
\end{eqnarray}
Change the dynamics of the VRRW in the
following way [recall (\ref{Etrans})]:
\begin{eqnarray}\label{Etransnew}
&&\mathbb{P}\bigl(X(t+1) =w |\mathcal{F}_t\bigr)\nonumber
\\[-8pt]\\[-8pt]
&&\qquad= \frac{Z(T_\varepsilon\wedge
t,w)}{\sum_{y\in
\{1,\ldots,d,\ell_1^i,\ldots,\ell_{r_i}^i\}:y \sim v}
Z(T_\varepsilon
\wedge t,y)} 1_{\{w \in
\{1,\ldots,d,\ell_1^i,\ldots,\ell_{r_i}^i\}\}}.\nonumber
\end{eqnarray}
In words, after
time $T_\varepsilon$ the step distribution does not anymore change
dynamically with the evolution of the walk; instead it is
``frozen'' to the configuration
\[
(Z(T_\varepsilon,1),\ldots,Z(T_\varepsilon,d),Z(T_\varepsilon
,\ell_1^1),\ldots,Z(T_\varepsilon,\ell_{r_d}^d)),
\]
and additional visits to the leaves attached at $g$ where $g\neq i$
become impossible.
Let
\[
\sigma_\varepsilon:=\inf\bigl\{ k\geq k_1 \dvtx  T_\varepsilon\leq\tau
_k^{(i)}\bigr\},
\]
and assume that we are given a family $\{U_k,k\geq k_1\}$ of independent
uniform $[0,1]$ random variables, and independent of the evolution of
the VRRW above.
Then define a modification $(\widetilde{X}'_k,\widetilde{Y}'_k)$ of
$(X',Y')$ by
$(\widetilde{X}'_{k_1},\widetilde{Y}'_{k_1})=(X'_{k_1},Y'_{k_1})$ and
%
\begin{equation}
\label{Etransnewad}
\hspace*{20pt}(\Delta\widetilde{X}'_k,\Delta\widetilde{Y}'_k):= \cases{
(\Delta X'_k,\Delta Y'_k),&\quad$k<\sigma_\varepsilon$,\cr
(d-1, 0),\qquad U_k \leq\widetilde{X}'_k/(\widetilde{X}'_k+\widetilde
{Y}'_k),&\quad $k\geq\sigma_\varepsilon$,\vspace*{2pt}\cr
(0, 1),\qquad U_k > \widetilde{X}'_k/(\widetilde{X}'_k+\widetilde
{Y}'_k),&\quad$k\geq\sigma_\varepsilon$.
}
\end{equation}
In words, the evolution of $(\widetilde{X}',\widetilde{Y}')$ is
identical to that of $(X',Y')$ up to time $\sigma_\varepsilon$, while
$(\widetilde{X}',\widetilde{Y}')$ evolves as the urn from Theorem
\ref{thurn1} from time $\sigma_\varepsilon$ onwards. In particular, the
asymptotic behavior of $(X',Y')$ and
$(\widetilde{X}',\widetilde{Y}')$ is the same on
$\{T_\varepsilon=\infty\}=\bigcap_{l\geq K} D_l'(\varepsilon)
\subset
\{\sigma_\varepsilon=\infty\}$.

The point of the above construction is that $(\widetilde
{X}',\widetilde
{Y}')$ satisfies
the hypotheses of
\cite{PV1999}, Lemma 3.5, with
\begin{eqnarray*}
a&=&1,\qquad b=b(\varepsilon)=\frac{d-1 + 3\varepsilon
d(d-3)/2}{1-3\varepsilon d/2}\quad \mbox{and }
\\
K&=&K(\varepsilon)= 2 \biggl(\frac{d-1 + 3\varepsilon
d(d-3)/2}{1-3\varepsilon d/2} \biggr)^2.
\end{eqnarray*}
Indeed, suppose $k<\sigma_\varepsilon$ (otherwise
the argument is trivial)
and note that then with probability
$Y'_k/(X'_k+Y'_k)=\widetilde{Y}'_k/(\widetilde{X}'_k+\widetilde{Y}'_k)$
we have
$X(\tau_k^{(i)}+1)\in\{\ell_1^i,\ldots,\ell_{r_i}^i\}$, so that
$(\Delta\widetilde{X}'_k, \Delta\widetilde{Y}'_k)=(0,1)$, while
with the remaining probability
$(\Delta\widetilde{X}'_k, \Delta\widetilde{Y}'_k)=(W_k,0)$
where $P(W_k\geq1)=1$ and conditionally on $\mathcal{F}_{\tau
_k^{(i)}}$, $W_k$ is
stochastically bounded from above by a Geometric random variable with
success probability
$(1-3\varepsilon d/2)/(d-1 + 3\varepsilon d(d-3)/2)$.
Here we use the
definition of the modified dynamics (\ref{Etransnew}) and (\ref
{Etransnewad}) together with
the fact (\ref{Eintbd}).

Due to \cite{PV1999}, Lemma 3.5, $(\widetilde{X}'_k/(\widetilde
{Y}'_k)^{b'}, k\geq k_1)$
is a positive supermartingale for any $b'>b(\varepsilon)$, hence converging,
and its limit
must be $0$, almost surely (strictly speaking, the supermartingale
property holds once
$\widetilde{Y}'_{k_1}$ is larger than some fixed
constant, but this we can assume WLOG).
Note that for any $\delta$ one can choose $\varepsilon>0$
sufficiently small so
that $d-1+\delta>b(\varepsilon)$.
Since $X'_\cdot/(Y'_\cdot)^{b'}$ and $\widetilde{X}'_\cdot
/(\widetilde
{Y}'_\cdot)^{b'}$
behave identically on $\{T_\varepsilon=\infty\}=\bigcap_{l\geq K}
D_l'(\varepsilon)$,
the statement of the lemma follows immediately from (\ref{EPDprime}).
\end{pf}

Now suppose that $\sum_{i=1}^d r_i>0$,
and denote by $\theta_k:=\sum_{i=1}^d L(t_k,i)/t_k>0$ the total
(rescaled) weight of the leaves.
Due to Lemma~\ref{lemrefinenew}, we have\break $\sum_{i=1}^d L(t_k,i)=o
(t_k^{1/(d-1)+\delta} )$,
hence
\begin{eqnarray}
\label{Eoneminusal}
\sum_{j=1}^d \alpha_j^{(k)}=1-\theta_k
\qquad\mbox{where }\theta_k=o \bigl(k^{-m [\vafrac{d-2}{d-1}-\delta]} \bigr).
\end{eqnarray}
Moreover, due to Lemma~\ref{lemrefinenew}, we have
$t_k^{1/(d-1)-\delta}=
o(\sum_{i=1}^d L(t_k,i))$,
therefore $t_k^{-(d-2)/(d-1)-\delta}=o(\theta_k)$
\begin{eqnarray*}
\Vert\pi(t)-\pi_{\mathsf{unif}}\Vert
&\geq&\sum_{i=1}^d \bigg|\frac
{Z(t_k,i)}{t_k}-\frac
{1}{d} \bigg| \geq\bigg|\sum_{i=1}^d \frac{Z(t_k,i)}{t_k} -1 \bigg|
\\
&=& \sum_{i=1}^d
\frac{L(t_k,i)}{t_k}
\gg t_k^{-\vafrac{d-2}{d-1}-\delta},\qquad \mbox{as }k\to\infty,
\end{eqnarray*}
yielding the lower bound claim (\ref{Eth2lower}) in Theorem \ref{th2}.

We continue toward the proof of (\ref{Eth2uppera}) and (\ref{Eth2upperb}).
Set
\[
\eta_k:= \eta(t_k)=1- d \min_{j=1,\ldots,d}\alpha_j^{(k)}\ge0,
\]
and let
\begin{eqnarray}\label{eqdeftbeta}
\tilde\beta=\min\biggl\{\beta,1,m \biggl(\frac{d-2}{d-1}-\delta\biggr)\biggr\},
\end{eqnarray}
where $\delta>0$ is very small.

The following statement is a
refinement of (3.28) in~\cite{V2001}.
\begin{lemma}\label{lemAnd}
On the event $E_k$ defined by (\ref{eqEkdef}) we have
\begin{eqnarray}\label{eqetareq}
\eta_{k+1}=\eta_k \biggl(1-\frac{mr(1-\eta_k) }{k} \biggr)+
\Theta\biggl(\frac1{k^{1+\tilde\beta}} \biggr),
\end{eqnarray}
where $r=r(k,\alpha^{(k)})\in[1/(d-1),1/(1-\eta_k)]$.
\end{lemma}
\begin{pf} Due to (\ref{Eoneminusal}) we have
\[
1-\sum_{j=1}^d \alpha_j^2\le1-\frac{ (\sum_{j=1}^d
\alpha_j )^2}d\le\biggl(1-\frac1d \biggr) +\frac{2\theta_k}d.
\]
Moreover, Theorem \ref{th1} implies in particular that $\mathbb
{P}(\bigcap
_{k\geq
k_0} \{\max_{i=1}^d \alpha_i^{(k)}<1/2\}) \to1$
as $k_0\to\infty$ (recall that $d\geq3$).
Since
$x \mapsto x(1-x)$ is an increasing function on $[0,1/2]$, we conclude that
asymptotically
\begin{eqnarray*}
1- \sum_{j=1}^d \alpha_j^2&=&\sum_{j=1}^d
\alpha_j (1-\alpha_j )+\theta_k\ge d\times\frac{1-\eta_k}d
\biggl(1-\frac{1-\eta_k}d \biggr)+\theta_k
\\&=& \biggl(1-\frac1d \biggr) - \biggl(1-\frac{2-\eta_k}d
\biggr)\eta_k+\theta_k.
\end{eqnarray*}
Thus we have shown
\begin{eqnarray}\label{eqcbound}
1-\sum_{j=1}^d \alpha_j^2&=& \biggl(1- \frac1d \biggr) - \frac{d-2}{d} \gamma
\eta_k+o\biggl(\frac
1{k^{m (\vafrac{d-2}{d-1}-\delta)}} \biggr)\nonumber
\\[-8pt]\\[-8pt]
\eqntext{ \mbox{where $\gamma\in [0,1+\eta_k/(d-2)]$}.\nonumber}
\end{eqnarray}
Note that
\begin{eqnarray*}
\alpha_i^{(k+1)}&=&\frac{\alpha_i k^m + N_k\afrac{\alpha_i
(1-\alpha_i)}{1-\sum_j \alpha_j^2}
+ O (k^{\vfrac{m-1}2+\nu} )
} {(k+1)^{m}},
\end{eqnarray*}
where the $O(\cdot)$ term comes from the estimation of the
event~(\ref{eqEkdef}). Thus
\begin{eqnarray*}
\alpha_i^{(k+1)}&=&\alpha_i \biggl[1-\frac mk +\frac mk \frac{1-\alpha
_i}{1-\sum
\alpha_j^2} \biggr]+O \biggl(\frac1{k^{1+\beta}} \biggr)
+\Theta\biggl(\frac1{k^2} \biggr)
\\
&=&\alpha_i \biggl[1+\frac mk \biggl(\frac{1-\alpha_i}{1-\sum
\alpha_j^2}-1 \biggr) \biggr]+\Theta\biggl(\frac1{k^{1+\tilde\beta}} \biggr).
\end{eqnarray*}
Since the last expression (without the $\Theta$ part) is
increasing in $\alpha_i$ for all sufficiently large $k$, it implies
that if $\alpha_i^{(k)}=\min_{j=1}^d\alpha_j^{(k)}$,
then $\alpha_i^{(k+1)}$ will again equal $\min_{j=1}^d\alpha_j^{(k+1)}$,
unless it
is ``overtaken'' by $\alpha_j^{(k+1)}$ for some other index $j$.
The latter case can happen only
if the difference $|\alpha_j^{(k)}- \alpha_i^{(k)}|$ is itself $O
(\frac
1{k^{1+\tilde\beta}} )$.
Hence it is always true that
\[
\min_{j=1}^d \alpha_j^{(k+1)}=\min_{i=1}^d \alpha_i^{(k)}
\biggl[1+\frac mk
\biggl(\frac{1-(\min_{i=1}^d \alpha_i^{(k)})}{1-\sum
\alpha_j^2}-1 \biggr) \biggr]+O \biggl(\frac
1{k^{1+\tilde\beta}} \biggr).
\]
This yields in turn
\begin{eqnarray*}
\eta_{k+1}&=&1-d \biggl(\frac{1-\eta_k}d \biggl[1+\frac mk \biggl(\frac{1-\vfrac{1-\eta_k}d}{1-\sum
\alpha_j^2}-1 \biggr) \biggr]+O \biggl(\frac
1{k^{1+\tilde\beta}} \biggr) \biggr)\nonumber
\\
&=&1-(1-\eta_k) \biggl[1+\frac mk \biggl(\frac{d-1 +\eta_k}{d-1-(d-2) \gamma
\eta
_k}-1 \biggr) \biggr]+O
\biggl(\frac
1{k^{1+\tilde\beta}} \biggr)
\\
&=&\eta_k \biggl(1- \frac{m(1-\eta_k)}k \times\frac{ 1+\gamma(d-2)
}{d-1-(d-2)\gamma\eta_k } \biggr)
+
O \biggl(\frac1{k^{1+\tilde\beta}} \biggr), \nonumber
\end{eqnarray*}
where
for the second equality we used (\ref{eqcbound}).
Since
\begin{eqnarray*}
\frac{d-1+\eta_k}{d-1 - (d-2)\eta_k-\eta_k^2}< \frac{1}{1-\eta_k},
\end{eqnarray*}
we get
\begin{eqnarray*}
\eta_{k+1}=\eta_k \biggl(1- \frac{m(1-\eta_k)r}k \biggr) +
O \biggl(\frac1{k^{1+\tilde\beta}} \biggr),
\end{eqnarray*}
where
$1/(d-1)\le r \le(1-\eta_k)^{-1}$.
\end{pf}

Recalling once again fact (\ref{Econvergence}) we can assume that for
$\varepsilon=1-2/d>0$ we have $\eta_k\le1-\varepsilon$, for all
large $k$.
This enables us applying Lemma~\ref{lemmaclaim1} with $C=m /(d-1)$.
Note that to get the best estimate of the speed of convergence we need
to make
$p(d,m):=\min\{C,\tilde\beta\}/m$ as large as possible, since
\[
\limsup_{k\to\infty} \eta_k h(k)=\limsup_{k\to\infty} \eta(k^m)
h(k)<\infty
\]
for an increasing function $h(\cdot)$ a.s.~implies
\[
\limsup_{t\to\infty} \eta(t) h(t^{1/m})<\infty.
\]
On the other hand, recalling the definition of $\tilde\beta$ from
(\ref{eqdeftbeta}), we have
\begin{eqnarray*}
p(d,m)&=&\min\biggl\{\frac1{d-1}, \frac{1}m, \frac{\beta}m, \frac
{d-2}{d-1}-\delta
\biggr\}
\\
&=&
\min\biggl\{\frac1{d-1}, \frac{1}m, \frac{1}2 - \frac{\delta_1+1/2}m,
\frac
{d-2}{d-1}-\delta\biggr\}.
\end{eqnarray*}
We can make $\beta$ as close as possible to $(m-1)/2$ by recalling
$\beta=(m-1)/2-\delta_1$, and taking $\delta_1>0$ arbitrarily small.
Similarly, $\delta>0$ can be made very small.
Given a particular choice of $\delta,\delta_1$,
observe that $\max_{m>1} p(d,m)$
is achieved at $3+2\delta_1$, so by setting $m=3+2\delta_1$ we obtain
\begin{eqnarray*}
p(d)&:=&p(d,3+2\delta_1)= \min\biggl\{\frac1{d-1}, \frac{1}{3+2\delta
_1}, \frac
{d-2}{d-1}-\delta\biggr\}
\\
&=&\min\biggl\{\frac1{d-1}, \frac{1}{3+2\delta_1},\frac1{d-1} + \biggl[ \frac
{d-3}{d-1}-\delta\biggr] \biggr\}
\\
&=&
\min\biggl\{\frac1{d-1}, \frac{1}{3+2\delta_1} \biggr\}.
\end{eqnarray*}
Consequently, $p(d)$ can be taken arbitrarily
close to $1/3$ if $d\in\{3,4\}$, while $p(d)=1/(d-1)$ for $d\ge5$.
Setting $C=3/(d-1)$ yields $\tilde\beta=\min\{1-\delta_1,1\}<C$ if
$d\in\{3,4\}$ and $\tilde\beta>C$ if $d\ge5$.
As already argued, this implies\break $\limsup\eta(t) t^{p(d)}<\infty$ due
to Lemma~\ref{lemmaclaim1},
and completes the proof of Theorem \ref{th2}.
\end{pf*}
\begin{rema}
There is a gap in the power between the upper and lower bounds on speed
of convergence in Theorem \ref{th2}.
One might wish to obtain further information on the lower bound using
(\ref{eqetareq}).
In fact, we would be able to conclude something provided
\begin{eqnarray*}
\eta_{k+1}\geq\eta_k \biggl(1-\frac{C(1-\eta_k) }{k} \biggr)+
\frac D{k^{1+\tilde\beta}},
\end{eqnarray*}
where both $C$ and $D$ are positive (or for $D$ negative,
under more complicated constraints on $C>0$ and $\tilde\beta$ that seem
difficult to verify).
Therefore, it is the lack
of knowledge of the sign (and magnitude) of the error term in (\ref
{eqetareq}) that
obstructs generalizing the above argument to obtaining lower bound estimate.
\end{rema}

\section*{Acknowledgments}

We are grateful to Robin Pemantle and to Pierre Tarr\`es for useful
discussions and for pointers to the literature.

%

\printaddresses

\end{document}